\documentclass[11pt,a4paper]{amsart}
\usepackage{amsmath,amsthm, amssymb}
\usepackage[text=16cm,centering]{geometry}
\usepackage{hyperref}
\usepackage{latexsym}\usepackage{color}
\usepackage{setspace, color}
\usepackage{bbm}
\usepackage{stmaryrd}

\newcommand{\NN}{\mathbb{N}}  \newcommand{\ZZ}{\mathbb{Z}}    
\newcommand{\RR}{\mathbb{R}}    \newcommand{\LL}{\mathbb{L}}

 \newcommand{\mcC}{\mathcal{C}}  \newcommand{\mcF}{\mathcal{F}}

        \newcommand{\mcA}{\mathcal{A}}
\newcommand{\mcE}{\mathcal{E}}

    \newcommand{\bfX}{{\bf X}}     \newcommand{\bfB}{\bf B} 
   \newcommand{\bfA}{\bf A} 

\newcommand{\bbX}{\mathbb{X}}   

\newcommand{\varep}{\varepsilon}
\newcommand{\ep}{\epsilon}
\renewcommand{\leq}{\leqslant}
\renewcommand{\geq}{\geqslant}

\theoremstyle{plain}
\newtheorem{theorem}{Theorem}
\newtheorem{lemma}[theorem]{Lemma}
\newtheorem{corollary}[theorem]{Corollary}

\theoremstyle{definition}
\newtheorem{definition}[theorem]{Definition}

\theoremstyle{remark}
\newtheorem{remark}[theorem]{Remark}

\numberwithin{equation}{section} 

\title[Unbounded rough drivers]{Unbounded rough drivers}
\date{\today}
\author[I. Bailleul]{Ismael Bailleul}
\address{IRMAR, 263 Avenue du General Leclerc, 35042 RENNES, France}
\email{ismael.bailleul@univ-rennes1.fr}
\author[M. Gubinelli]{Massimiliano Gubinelli}
\address{Institut Universitaire de France}
\address{CEREMADE \& CNRS UMR 7534,  Universit\'e Paris Dauphine, Place du Mar\'echal De Lattre De Tassigny, 75775 PARIS cedex 16}
\email{massimiliano.gubinelli@ceremade.dauphine.fr}

\begin{document}

\begin{abstract}
We propose a theory of linear differential equations driven by unbounded operator-valued rough signals. As an application we consider rough linear transport equations and more general linear hyperbolic symmetric systems of equations driven by time-dependent vector fields which are only distributions in the time direction. 
\end{abstract}

\maketitle

\section{Introduction}
\label{SectionIntroduction}

In this paper we start the program of developing a general theory of rough PDEs aiming at extending classical PDE tools such as weak solutions, a priori estimates, compactness results, duality. This is a quite unexplored territory where few tools are available, so as a start, we will content ourselves in this work with the study of linear symmetric hyperbolic systems of the form
\begin{equation}
\label{eq:basic}
\partial_t f + a\nabla f = 0,
\end{equation}
where $f$ is an $\RR^N$-valued space-time distribution on $\RR_+\times\RR^d$, and $a :\RR_+\times \RR^d \to  \textrm{L}(\RR^d, \RR^{N\times N})$ is a $N\times N$ matrix-valued family of time--dependent vector fields in $\RR^d$. This setting includes as a particular case scalar transport equations. Moreover we restrict our attention to the case where the matrix-valued vector field $a$ is only a distribution in the time variable, rather than a regular bounded function. We however retain some smoothness assumption in the space variable, as expected from the fact that general transport equations do not possess the regularisation properties needed to drive them with space-time irregular signals. Even in the classical setting it is known that non-regular coefficients  can give rise to non-uniqueness of weak solutions~\cite{FlandoliMilano}.

\smallskip

When $a$ is only a distribution in the time variable the above weak formulation is not available since in the classical setting solutions are considered in spaces like $C([0,T],\LL^2(\RR^d))$ for general symmetric systems  or $L^\infty(\RR\times \RR^d)$ for scalar transport equation. In this case, the product $a\nabla f$ is not well-defined, not even in a distributional setting. Rough paths have their origin in the need to handle such difficulties in the case of ordinary differential equations driven by distribution valued signals \cite{Lyons98,Lyons2002,Lyons2007,Friz2010}. Controlled rough paths have been introduced in \cite{Gubinelli04} as a setting for considering more general problems; they were used successfully in the study of some stochastic partial differential equations~\cite{Bessaih2005,GT,Deya2012,Hairer2011,Gubinelli2012a,Hairer2012,Hairer2013},  including the remarkable solution by Hairer of the KPZ equation \cite{HairerKPZ}. 

\smallskip

These developments have ultimately lead to the notion of paracontrolled distributions introduced by Gubinelli, Imkeller and Perkowski~\cite{GIP} and to Hairer's general theory of regularity structures~\cite{HairerRegularity},  providing a  framework for the analysis of non-linear operations on distributions. Despite their successes, these new tools and methods are somehow designed to deal with a prescribed class of singular PDEs which is far from exhausting the set of all interesting singular PDEs. 

\smallskip

PDEs with irregular signals have been studied using directly rough path methods also by Friz and co-authors~\cite{Caruana2009, CaruanaFrizOberhauser, FrizOberhauser1, FrizOberhauser2, FrizGess}. They have developed an approach to some  fully non-linear PDEs and conservation laws driven by rough signals by interpreting these equations  as transformations of classical PDEs, generalising the method of characteristics. Subsequently a combination of rough path stability and PDE stability allows to go from smooth signal to the wider class of signals described by rough paths. Entropy solutions to scalar conservation laws with rough drivers have also been analysed also by P.-L. Lions, Souganidis and coauthors~\cite{LionsPerthameSouganidis1,LionsPerthameSouganidis2, GessSouganidis}. A major drawback of this otherwise effective approach is that there is no intrinsic notion of solution to the PDE and that the study of the properties of the PDE has to be done on a global level.

\smallskip

In recent works intrinsic notions of weak solution to rough PDEs have been proposed by Tindel, Gubinelli and Torecilla~\cite{GTT} for viscosity solutions,  Catellier~\cite{CatellierTransport} for weak solutions to linear transport equations (see also Hu and Le~\cite{HuLe} for classical solutions to transport equations) and more recently by Diehl, Friz and Stannat~\cite{DiehlFrizStannat} for a general class of parabolic equations. All these notions are based on a weak formulation of the equation where the irregularity of some data is taken into account via the framework of controlled paths introduced in~\cite{Gubinelli04}. However in all these papers explicit formulas involving the flow of rough characteristics play an important role, and this sets apart the study of rough PDEs with respect to the study of weak solutions to more regular PDEs. One of the main motivations of our investigations is an effort to understanding what kind of robust arguments can be put forward in the study of the a priori regularity of distributions satisfying certain rough PDEs formulated in the language of controlled paths. Extensions to regularity structures or paracontrolled distributions will be considered in forthcoming work. 

\bigskip

We study equation \eqref{eq:basic} by working in the technically easier setting of controlled paths. To motivate our formalism, note that a formal integration of the weak formulation \eqref{eq:basic} over any time interval $[s,t]$, gives an equation of the form
$$
f_t = f_s + \int_s^t V_r f_r dr, 
$$ 
where $V_r = a_r\nabla$ is a matrix-valued vector-field and $f_r(x) = f(r,x)$, is a convenient notation of the distribution $f$ evaluated at time $r$, assuming this make sense. An expansion for the time evolution of $f$ is obtained by iterating the above equation, and reads
\begin{equation}
\label{eq:basic-increments}
f_t = f_s + A^1_{ts} f_s + A^2_{ts} f_s + R_{ts},
\end{equation}
where
$$
A^1_{ts} = \int_s^t V_r  dr,\qquad\textrm{ and }\qquad A^2_{ts} = \int_s^t \int_s^r V_r V_{r'} dr' dr,
$$
are respectively a first order differential operator (that is a vector field) and a second order differential operator, for each $s\leq t$. As a function of $(s,t)$, they satisfy formally \textit{Chen's relation}
\begin{equation}
\label{eq:operator-chen}
A^2_{ts} =A^2_{tu}+ A^2_{us} + A^{1}_{tu}A^{1}_{us}
\end{equation}
for all $0\leq s \leq u \leq t$. It is a key observation of rough path theory that equation \eqref{eq:basic-increments} can be used as a replacement for the differential or integral formulation 
of equation \eqref{eq:basic} if the remainder term can be shown to be sufficiently small as $t-s$ goes to $0$, in a suitable sense.

\medskip

We shall call a \emph{rough driver} an operator-valued $2$-index maps ${\bfA}_{ts}=(A^1_{ts},A^2_{ts})$ satisfying the operator Chen relation \eqref{eq:operator-chen} and some regularity assumptions. Building on the above picture, a path with values in some Banach space where the operators $A^1_{ts}$ and $A^2_{ts}$ act, will be said to solve the rough linear equation
$$
df_s = {\bfA}(ds)\,f_s.
$$
if the Taylor expansion \eqref{eq:basic-increments} holds.

\smallskip

There is a complete theory of such equations in the case where the equation is set in a Banach algebra and the operators $A^1_{ts}$ and $A^2_{ts}$ are given by left multiplication by some elements of the algebra. It is however natural, in the present PDE setting, to consider also unbounded operators $A^1, A^2$, which makes the use of rough paths ideas non-trivial, unless we work in the analytic category or in similar topologies. 

\medskip
 
We lay out in this work a theory of such rough linear equation driven by unbounded drivers $\bfA$, and obtain some a priori estimates that are used to study the well--posedness of some classes of linear symmetric systems in $\LL^2$ and of the rough transport equation in $\LL^\infty$. The major difficulty which has to be overcome is the lack of a Gronwall lemma for rough equations and the main contribution to this paper is to develop suitable a priori estimates on weak controlled solutions that replace the use of Gronwall lemma in classical proofs. Along the way we refine the standard theory of controlled path by introducing weighed norms compatible with the sewing map and by revisiting the theory of linear rough differential equations in the context of bounded drivers. 

\medskip

As a guide for the reader, here is how we have organised our work. Section~\ref{SectionWeightedNorms} provides a refined version of the sewing lemma that allows to keep track of the growth in time of the additive function associated with an almost-additive $2$-index map. This result is used in Section~\ref{SectionBoundedDrivers} in the proof of the well-posed character of linear differential equations driven by the bounded rough drivers defined there. \emph{Unbounded rough drivers} are introduced in Section~\ref{SectionUnboundedRoughDrivers}, where some fundamental a priori estimate is proved. An $\LL^2$ theory of rough linear equations is developed for a class of unbounded drivers, that contains as a particular example the rough linear transport equation. Our main workhorse here is a novel renormalisation lemma obtained from a tensorization argument in the line of the "doubling of variables" method commonly used in the setting of transport equations or conservation laws. A complete $\LL^\infty$ theory of rough transport equations is given in Section~\ref{SectionLInftyTheory}.

\bigskip

\noindent \textbf{Acknowledgments} -- The authors would like to express their gratitude to Martina Hofmanova and Mario Maurelli who discovered a major error in the first version of the paper and hinted to a possible strategy to overcome it.
\bigskip

\noindent \textbf{Notations} -- We gather here for reference a number of notations that will be used throughout the text.
\begin{itemize}
   \item We shall denote by E a generic Banach space. Given two Banach spaces E and F, we denote by $\textrm{L}(\textrm{E},\textrm{F})$ the set of continuous linear maps from E to F. \vspace{0.1cm}
   \item We shall denote by $c$ a constant whose value may change from place to place. \vspace{0.1cm}
   \item Given two positive real numbers $\alpha,\beta$, we shall write $\alpha \lesssim \beta$ to say that $\alpha\leq c\beta$, for some positive constant $c$. To indicate that this constant $c$ depends on a parameter $\lambda$, we write $\alpha \lesssim_\lambda \beta$. \vspace{0.1cm}
   \item Denote by $\|\cdot\|_{\alpha\,;\,k}$ the $\alpha$-H\"older norm of an $E_k$-valued path, for $k\in\ZZ$, and by $\|\cdot\|_{\alpha\,;\,\textrm{E}}$ the $\alpha$-H\"older norm of an E-valued path.
\end{itemize}

\bigskip

\section{Weighted norms}
\label{SectionWeightedNorms}

We introduce  some weighted norms that will be useful in getting a priori estimates on the growth of solutions to the linear differential equations studied in Section~\ref{SectionBoundedDrivers}. These norms are modelled on Picard's well-known norms
$$
\llparenthesis\, f \,\rrparenthesis := \underset{t\geq 0}{\sup}\,e^{-\lambda^{-1} t}| f(t) |,
$$
introduced in the study of ordinary differential equations in order to provide a functional setting where to get well--posedness results on the whole time interval $[0,\infty)$, as a direct consequence of the Banach fixed point theorem, and to get as a consequence a control on the growth of the size of solutions. 

\medskip

Let $T$ be a possibly infinite positive time horizon. As is common in rough paths theory, we shall work with Banach space-valued multi-index maps, mainly $2$ and $3$-index maps, defined on the simplexes
$$
\{(s,t)\in [0,T)^2\,;\,s\leq t\}\quad \textrm{ and } \quad \{(s,u,t)\in [0,T)^3\,;\,s\leq u\leq t\}.
$$
With Picard's norm in mind, we introduce a norm on the set of $2$ and $3$-index maps which captures both their H\"older size and their growth at infinity. 
Given $\lambda >0$, an increasing non-negative function $g$ defined on $\RR_+$ and a non-negative H\"older exponent $\gamma$, we define the $(\gamma,g)$-norm of a $2$-index map $a$, and a $3$-index map $b$, by the formulae
$$
\llparenthesis\, a \,\rrparenthesis_{\gamma,g} := \underset{\stackrel{0\leq s<t<T}{|t-s|\le \lambda}}{\sup}\;\,\frac{| a_{ts}|}{g(\lambda^{-1} t)\,|t-s|^\gamma},
$$
and 
$$
\llparenthesis\, b \,\rrparenthesis_{\gamma,g} := \underset{\stackrel{0\leq s< u<t<T}{{|t-s|\le \lambda}}}{\sup}\;\,\frac{| b_{tus}|}{g(\lambda^{-1} t)\,|t-s|^\gamma};
$$
note the following comparison: For $0\leq \gamma'\leq \gamma$, we have
\begin{equation}
\label{EqComparisonNorms}
\llparenthesis\, \cdot \,\rrparenthesis_{\gamma',g} \leq \lambda^{\gamma-\gamma'} \llparenthesis\, \cdot \,\rrparenthesis_{\gamma,g}.
\end{equation}

\smallskip

Recall that given an E-valued $2$-index map $a$, the \textit{sewing lemma} \cite{Gubinelli04, FdlP} asserts that if the inequality
\begin{equation}
\label{EqAlmostAdditivity}
| a_{ts} - (a_{tu}+a_{us})| \leq c|t-s|^\zeta
\end{equation}
holds for all $0\leq s\leq u\leq t<T$, with $t-s\leq 1$ say, for some exponent $\zeta>1$ and some positive constant $c$, then there exists a unique map $A : [0,T)\rightarrow \textrm{E}$ whose increments $\delta A_{ts} := A_t-A_s$ are well-approximated by $a_{ts}$, in the sense that 
$$
|\delta A_{ts} - a_{ts} | \lesssim |t-s|^\zeta,
$$
for all $t-s\leq 1$ say. Moreover, it $t_i$ denotes the times of a finite partition $\pi_{ts}$ of an interval $(s,t)$, with mesh $|\pi_{ts}|$, we have
\begin{equation}
\label{EqEstimateAlmostAdditiveFunctional}
|\delta A_{ts} - \sum a_{t_{i+1}t_i}| \lesssim |t-s||\pi_{ts}|^{\zeta-1}.
\end{equation}

\smallskip

The \textbf{sewing map} associates to the above $2$-index map $a$ the $2$-index map
$$
\Lambda(a)_{ts} := \delta A_{ts} - a_{ts}.
$$
For a given function $g:\RR\to\RR$ let 
$$
G(t) = g(t) + \int_0^t g(r)\,dr,
$$
and write as usual  $\delta a_{tus}$ for $a_{ts} - (a_{tu}+a_{us})$, for any $0\leq s\leq u\leq t<T$. The following lemma provides an estimate of the weighted norm of $\Lambda(a)$ in terms of the weighted norm of $\delta a$; an estimate for the growth of $A$ follows as a consequence. 

\begin{lemma}
\label{LemmaSewingMap}
There exists a positive constant $c_\zeta$, depending only on $\zeta$, such that 
$$
\llparenthesis\,\Lambda(a)\,\rrparenthesis_{\zeta,g} \leq c_\zeta \, \llparenthesis\,\delta a\,\rrparenthesis_{\zeta,g}.
$$
\end{lemma}

\medskip

\begin{proof}
Given $0\leq s<t$ and $n\geq 1$, set $t_i = s+i2^{-n}(t-s)$, and note that the telescopic sum
\begin{equation*}
(\star) := \sum_{i=0}^{2^n-1}a_{t_{i+1}t_i}\; - a_{ts} = \sum_{k=0}^{n-1}\;\sum_{\ell=0}^{2^{n-(k+1)}} (a_{(\ell+2)2^k\,(\ell+1)2^k}+a_{(\ell+1)2^k\,\ell2^k} - a_{(\ell+2)2^k\,\ell2^k} )
\end{equation*}
provides a control of the quantity $(\star)$ in terms of $\delta a$ only. Sending $n$ to infinity, we see by identity \eqref{EqEstimateAlmostAdditiveFunctional} that there exists a constant $c_\zeta$ depending only on $\zeta$, such that 
$$
| \Lambda(a)_{ts} | \leq c_\zeta \, |t-s|^\zeta \,\sup\left\{\frac{|\delta a_{t'u's'}|}{|t'-s'|^\zeta}\,;\,0\leq s<s'<u'<t'\leq t\right\}.
$$
If $|t-s| \leq \lambda$ the above inequality implies clearly that we have
\begin{equation}
\label{EqEstimateLambdaA}
\frac{| \Lambda(a)_{ts} |}{g(\lambda^{-1} t)\,|t-s|^\zeta} \leq c_\zeta \, |\delta a|_{\zeta,g}.
\end{equation}
\end{proof}

\medskip

We will consider only the particular choice of function $g(t) = e^{t}$ and we shall set for a path $f : [0,T)\rightarrow \textrm{E}$
\begin{equation}
\label{EqDefnNormFunction}
\llparenthesis\, f_\bullet \,\rrparenthesis :=  \underset{t\geq 0}{\sup}\,e^{-\lambda^{-1} t}\,| f(t) |,
\end{equation}
and for a $2$-index map $a$, and a  $3$-index map $b$,
\begin{equation}
\label{EqDefnNormMultiIndexMaps}
\llparenthesis\, a \,\rrparenthesis_\gamma := \underset{\stackrel{0\leq s<t<T}{|t-s|\le \lambda}}{\sup} \;\, e^{-\lambda^{-1} t} \,\frac{| a_{ts}|}{|t-s|^\gamma}\,, \quad \textrm{ and }\quad 
\llparenthesis\, b \,\rrparenthesis_\gamma := \underset{\stackrel{0\leq s< u<t<T}{{|t-s|\le \lambda}}}{\sup} \;\, e^{-\lambda^{-1} t} \,\frac{| b_{tus}|}{|t-s|^\gamma}\,.
\end{equation}
Note that these norms depend on a choice of parameter $\lambda$, which may be tuned on demand. This will be particularly useful in the statement and proof of Theorem~\ref{ThmIntegrationLinearBounded} below, giving the well-posed character of some linear rough differential equation. Note also that we can compare $\llparenthesis\, f_\bullet \,\rrparenthesis$ and $\llparenthesis\, \delta f \,\rrparenthesis_\gamma$, as the inequality
\begin{equation}
\label{EqUsefulEstimate}
\llparenthesis\, f_\bullet \,\rrparenthesis \leq |f_0| + e \lambda^\gamma\,\llparenthesis\, \delta f \,\rrparenthesis_\gamma
\end{equation}
holds for all $\lambda>0$. The proof is easy. Given a fixed time $t$, let $\{0=t_0< t_1 < \cdots < t_n = t\}$ be a partition of the interval $[0,t]$ into sub-intervals size at most $\lambda$. Then 
$$
|f(t)| \le |f(0)| + \sum_{k=0}^{n-1} | f(t_k) - f(t_{k+1}) | \le  |f(0)| +  \lambda^\gamma  \llparenthesis\, \delta f \,\rrparenthesis_\gamma \sum_{k=0}^{n-1} e^{\frac{t_{k+1}}{\lambda}}.
$$
But now
$$
\sum_{k=0}^{n-1} e^{\frac{t_{k+1}}{\lambda}}
\le   \lambda^{-1}\int_{t_{k}+\tau}^{t_{k+1}+\lambda} e^\frac{s}{\lambda}\, ds \le  \lambda^{-1} \int_0^{t+\lambda} e^\frac{s}{\lambda}\, ds \le e^{\frac{s}{\lambda}+1}.
$$
The comparison estimate \eqref{EqUsefulEstimate} will be our starting point in the proof of the a priori estimate \eqref{EqLinearBound} in Theorem~\ref{ThmIntegrationLinearBounded} below. 

\smallskip

Last, a $2$-index map $a$ such that $\sup_{t-s\leq 1}\,\frac{| a_{ts} |}{|t-s|^\gamma}$ is finite will be called a $\gamma$\textbf{-H\"older map.}

\section{Linear differential equations with bounded rough drivers}
\label{SectionBoundedDrivers}

Let $(\mcA,|\cdot|)$ be a Banach algebra with unit ${\bf 1}_{\mathcal{A}}$; one may think for instance to the space of continuous linear maps from some Hilbert space to itself, or to the truncated tensor algebra over some Banach space, equipped with a tensor norm and completed for that norm. We introduce in this section a notion of bounded rough driver in the Banach algebra $\mcA$, and show that they generate some flows on the algebra. 

\begin{definition}
Let $\frac{1}{3}<\gamma\leq \frac{1}{2}$. A \textbf{\emph{bounded $\gamma$-rough driver in $\mcA$}} is a pair ${\bf A}= (A^1,A^2)$ of $\mcA$-valued $2$-index maps satisfying \emph{Chen's relations}
\begin{equation}
\label{EqChenRelations}
\delta A^1=0, \quad \textrm{ and }\quad \delta A^2_{tus} = A^1_{tu}A^1_{us},
\end{equation}
and such that $A^1$ is $\gamma$-H\"older and $A^2$ is $2\gamma$-H\"older. The norm of $\bf A$ is defined by the formula
$$
\| {\bf A} \| := \sup_{0\leq s<t<T\,;\,t-s\leq 1}\;\frac{|A^1_{ts}|}{|t-s|^\gamma}\vee\frac{|A^2_{ts}|}{|t-s|^{2\gamma}}.
$$
\end{definition}

\medskip

As an elementary example, think of $\mcA$ as the truncated tensor algebra $\oplus_{i=0}^N(\RR^\ell)^{\otimes i}$ over $\RR^\ell$, for $N\geq 2$, and consider a weak geometric $\gamma$-H\"older rough path ${\bfX}_{ts} = 1\oplus X_{ts}\oplus \bbX_{ts}\in \oplus_{i=0}^N(\RR^\ell)^{\otimes i}$, with $2\leq \gamma<3$. Left multiplication by $X_{ts}$ and $\bbX_{ts}$ define operators $A^1$ and $A^2$ that are the components of a rough driver.

\medskip

Recall the weighted norms $\llparenthesis\,\cdot\,\rrparenthesis$ and $\llparenthesis\,\cdot\,\rrparenthesis_\gamma$ defined by identities \eqref{EqDefnNormFunction} and \eqref{EqDefnNormMultiIndexMaps}, respectively, depend on some parameter $\lambda$. The proof of the following well--posedness result for linear rough differential equations driven by bounded $\gamma$-rough drivers shows the interest of being flexible on the tuning of $\lambda$.

\begin{theorem}[Integration of bounded rough drivers]
\label{ThmIntegrationLinearBounded}
Given any initial condition $J_0\in\mcA$, there  exists a unique $\gamma$-H\"older path $f_\bullet$ starting from $f_0$, and such that the formula 
  \begin{equation*}
    \delta f_{t, s} - (A^1_{t, s} +A^2_{t, s}) f_s 
  \end{equation*}
defines a $3\gamma$-H\"older $2$-index map $f^\natural$. Moreover, the following estimate holds
  \begin{equation}
  \label{EqLinearBound}
    \llparenthesis\, f_\bullet \,\rrparenthesis \leqslant 2\,| f_0 |
  \end{equation}
for all $\lambda$ greater than some $\lambda_0$ depending only on $\| \bf A \|$ and $\gamma$. When $f_0  ={\bf 1}_{\mathcal{A}}$, we will use the notation $ f_t =e^{\bf A}_{t,0}$; the flow property
$$
e^{\bf A}_{t, s} = e^{\bf A}_{t, u} e^{\bf A}_{u, s},
$$
holds for all $0\leq s\leq u\leq t<T$.
\end{theorem}

\medskip

Applied to the above example of rough driver in the truncated tensor product space, this well--posedness result provides a proof of Lyons' extension theorem \cite{Lyons98}.

\medskip

\begin{proof}
{\bf a) A priori estimate --} Let us prove the a priori bound \eqref{EqLinearBound} first; the uniqueness claim in the theorem follows from this bound and the linear character of the problem. As mentioned above, we start from the inequality
$$
\llparenthesis\,f_\bullet \,\rrparenthesis \leq | f_0 | + c \lambda^\gamma \llparenthesis\, \delta f \,\rrparenthesis _\gamma,
$$
and try and write $\llparenthesis\, \delta f \,\rrparenthesis _\gamma$ in terms of $\llparenthesis\,f_\bullet \,\rrparenthesis$; this can be done as follows.

By using Chen's relation \eqref{EqChenRelations} and the definition of the remainder $f^\natural$, the identity
$$
- \delta f^\natural_{t, u, s} = A^1_{t, u} (\delta f_{u, s}  -A^1_{u, s} f_s) +A^2_{t, u} \delta f_{u, s} = A^1_{t, u} (A^2_{t, s} f_s + f^{\natural}_{s, t}) + A^2_{t, u} \delta f_{u, s}, 
$$
gives us the estimate
$$
\llparenthesis\, \delta f^{\natural} \,\rrparenthesis_{3 \gamma} \leq \| {\bf A} \| ^2 \llparenthesis\, f_\bullet \,\rrparenthesis + \| {\bf A} \| \llparenthesis\, f^{\natural} \,\rrparenthesis_{2 \gamma} + \| {\bf A} \| \llparenthesis\, \delta f \,\rrparenthesis_\gamma.
$$ 
But since $3 \gamma > 1$, we have 
$$
f^{\natural} = \Lambda \delta f^{\natural},
$$ 
so the inequality
$$
\llparenthesis\, f^{\natural} \,\rrparenthesis_{3 \gamma} \lesssim_\gamma \| {\bf A} \|^2 \llparenthesis\, f_\bullet \,\rrparenthesis + \| {\bf A} \| \llparenthesis\, f^{\natural} \,\rrparenthesis_{2 \gamma} + \| {\bf A} \| \llparenthesis\, \delta f \,\rrparenthesis_\gamma
$$
follows from Lemma~\ref{LemmaSewingMap}. Using the inequality $\llparenthesis\,  f^{\natural} \,\rrparenthesis_{2 \gamma} \leq \lambda^\gamma \llparenthesis\,  f^{\natural} \,\rrparenthesis_{3 \gamma}$, emphasized in \eqref{EqComparisonNorms}, the above equation gives 
$$
\llparenthesis\,  f^{\natural} \,\rrparenthesis_{3 \gamma} \lesssim_\gamma \| {\bf A} \|^2\, \llparenthesis\,  f_\bullet \,\rrparenthesis  + \lambda^\gamma \| {\bf A} \|\,\llparenthesis\,  f^{\natural}\,\rrparenthesis_{3 \gamma} + \| {\bf A} \|\, \llparenthesis\, \delta f \,\rrparenthesis_\gamma.
$$
For $\lambda$ small enough so that $\lambda^\gamma \| {\bf A} \| \leq \frac{1}{2}$, we obtain
$$
\llparenthesis\,  f^{\natural} \,\rrparenthesis_{3 \gamma} \lesssim_{\gamma} \| {\bf A} \|^2\, \llparenthesis\,  f_\bullet \,\rrparenthesis  + \| {\bf A} \|\, \llparenthesis\, \delta f \,\rrparenthesis_\gamma,
$$
so, using again the definition of the remainder $f^\natural$, and the observation that 
$$
\llparenthesis\, A^2 f\,\rrparenthesis_{\gamma} \lesssim_{\gamma}  \lambda^\gamma \llparenthesis\, A^2 f\,\rrparenthesis_{2 \gamma} \lesssim_{\gamma} \lambda^\gamma \| {\bf A} \|\, \llparenthesis\, f_\bullet \,\rrparenthesis,
$$
we obtain the estimate
\begin{equation*}
\begin{split}
\llparenthesis\, \delta f \,\rrparenthesis_\gamma &\lesssim_{\gamma} \| {\bf A} \|\llparenthesis\,f_\bullet \,\rrparenthesis + \llparenthesis\, A^2 f\,\rrparenthesis_\gamma + \llparenthesis\, f^{\natural}\,\rrparenthesis_\gamma \\
                            &\lesssim_{\gamma} \{\| {\bf A} \| (1 + \lambda^{2 \gamma}\| {\bf A} \|) + \lambda^\gamma \| {\bf A} \|\} \llparenthesis\, f_\bullet \,\rrparenthesis + \lambda^{2 \gamma}\| {\bf A} \|\,\llparenthesis\, \delta f \,\rrparenthesis_{\gamma}.
\end{split}
\end{equation*}
Taking $\lambda$ small enough, depending only on $\| \bf A \|$, we eventually see that 
$$ 
\llparenthesis\,\delta f \,\rrparenthesis_\gamma \lesssim_{\gamma}\{\| {\bf A} \| (1 +  \lambda^{2 \gamma} \| {\bf A} \|) + \lambda^\gamma  \| {\bf A} \|\} \llparenthesis\, f_\bullet\,\rrparenthesis.
$$
The a priori estimate $\llparenthesis\, f _\bullet\,\rrparenthesis \leqslant 2 | f_0 |$, follows now from a choice of sufficiently small  parameter $\lambda$, since $\llparenthesis\,f_\bullet \,\rrparenthesis \leq | f_0 | + c \lambda^\gamma \llparenthesis\, \delta f \,\rrparenthesis _\gamma$.

\bigskip

{\bf b) Existence --} We can run a Picard iteration to prove the existence of a path satisfying the conditions of the theorem. Set first
$$
f^0_t = f_0,\quad \textrm{ and } \quad f^1_t =A^1_{t0} f_0
$$ 
for all $t \in \RR$. Given the paths $f^{n-1}_\bullet, f^n_\bullet$, the $2$-index map 
$$
a^n_{ts} := A^1_{ts} f^n_s +A^2_{ts} f_s^{n - 1}
$$
satisfies the almost-additivity condition \eqref{EqAlmostAdditivity} with $\zeta=3\gamma>1$ here, so there is, by the sewing lemma, a unique $\gamma$-H\"older path $f^{n + 1}_\bullet$ for which the formula
$$ 
\delta f^{n + 1}_{ts} - (A^1_{ts} f^n_s +A^2_{ts} f_s^{n - 1}) 
$$
defines $2$-index $3\gamma$-H\"older map $f^{n + 1, \natural}$. Setting $g^0_\bullet := 0$ and
$$
g^{n + 1}_\bullet := f^{n + 1}_\bullet - f^n_\bullet, \quad  g^{n + 1, \natural} := f^{n + 1, \natural} - f^{n, \natural},
$$ 
for all $n \geq 0$, we have
$$
\delta g^{n + 1}_{ts} =A^1_{ts} g^n_s +A^2_{ts} g_s^{n - 1} + g^{n + 1, \natural}_{ts}. 
$$
Note moreover that we have the identity
\begin{equation*}
\begin{split}
- \delta g^{n + 1, \natural}_{t, u, s} &= A^1_{t, u} (\delta  g^n_{u, s} -A^1_{u, s} g^{n - 1}_s) +A^2_{t, u} \delta g^{n - 1}_{u, s} \\
&= A^1_{t, u} (A^2_{t, s} g^{n - 2}_s + g^{n, \natural}_{s, t}) +A^2_{t, u} \delta g^{n - 1}_{u, s},
\end{split}
\end{equation*}
so, proceeding as in the proof of the a priori bound, we see that the inequality
$$
\llparenthesis\, g^{n - 1}_\bullet\,\rrparenthesis + \llparenthesis\, \delta g^n \,\rrparenthesis_\gamma + \llparenthesis\, g^{n + 1, \natural} \,\rrparenthesis_{3 \gamma} \lesssim_{\gamma, \| {\bf A} \|} \lambda^\gamma \{\llparenthesis\, g^{n - 2}_\bullet \,\rrparenthesis + \llparenthesis\, \delta g^{n - 1} \,\rrparenthesis_\gamma + \llparenthesis\, g^{n, \natural} \,\rrparenthesis_{3 \gamma}\}
$$  
holds, by choosing $\lambda$ small enough. The estimate
$$
\llparenthesis\, g^{n - 1}_\bullet\,\rrparenthesis + \llparenthesis\, \delta g^n \,\rrparenthesis_\gamma + \llparenthesis\, g^{n + 1, \natural} \,\rrparenthesis_{3 \gamma} \lesssim_{\gamma, \| {\bf A} \|} \lambda^{\gamma n}
$$
follows as a consequence, so the series $f^n = f^0 + \sum_{n \geqslant 1} g^n$ converges in the complete space of $\mcA$-valued $\gamma$-H\"older paths, and defines a path in $\mathcal{A}$ satisfying the conditions of the theorem.
\end{proof}

\smallskip

\begin{remark}
\begin{enumerate}
   \item Note that the proof given above gives back the known sharp growth rate $\exp((2\| {\bf A}\|)^\gamma\,t)$ for $|f_t|$; see \cite{FrizOberhauser}. Bounded rough drivers can also be integrated by defining recursively the $(n\gamma)$-H\"older $\mcA$-valued $2$-index map $A^n$ using the formula
$$
\delta A^n_{t, u, s} = \sum_{k = 1}^{n - 1} A^{n - k}_{t, u} A^k_{u, s},
$$
and setting 
$$
e^{\bf A}_{t, s} = \sum_{n = 0}^{\infty} A^n_{t, s}.
$$ 
Standard estimates on the sewing map \cite{Gubinelli04} show that $\delta A^n$ has $n \gamma$-H\"older norm  no greater than $(n!)^{- \gamma}$, so the above series converges in $\mcA$ for all $0 \leq s \leq t$. The flow property is obtained by a direct calculation, and setting $f_t := e^{\bf A}_{t, 0} f_0$, we see that the path $f_\bullet$ solves the problem.   \vspace{0.2cm}

   \item {\bf Linear rough differential equations with a linear drift -- } The above theory extends easily to rough equations of the form
\begin{equation}
\label{EqRoughLinearDrift}
\delta f_{t, s} = \int_s^t B_r f_r\,dr + A^1_{t, s} f_s + A^2_{t, s} f_s + f^{\natural}_{t, s}
\end{equation}
where $B \in \LL^{\infty} (\RR; \mcA)$ is a bounded measurable family of bounded operators. This equation is the rigorous meaning to give to solutions of the differential equation
$$
\frac{d}{dt} f = (B_r + \dot{A}_r) f_r
$$
where $\dot{A}_t = \partial_t A_{t0}^1$. In case $\dot{A}_t \in \LL^{\infty} (\RR; \mcA)$, the two formulations are equivalent provided ${\bf A} = (A^1, A^2)$ is defined by the formula
$$ 
A^1_{t, s} = \int_s^t \dot{A}_r \,dr, \quad A^2_{t, s} = \int_s^t \int_s^r \dot{A}_u \dot{A}_r \,du dr.
$$
The proof of Theorem~\ref{ThmIntegrationLinearBounded} can be easily adapted in the present setting, and the final lower bound on $\lambda$ gets an additional dependence on $\| B \|_\infty$. It yields moreover the following Duhamel formula
$$
f_t = e^{\bf A}_{t, 0} f_0 + \int_0^t e^{\bf A}_{t, r} B_r f_r\,dr.
$$
Indeed, let $f_\bullet$ be a function satifsying the above identity. If one computes the increment of the right hand side in the above equation, we get
$$
\delta f_{t, s} = \int_s^t e^{\bf A}_{t, r} B_r f_r \,dr - (e^{\bf A}_{t, s} -\textrm{\emph{Id}}) f_s = \int_s^t B_r f_r \,dr + A^1_{t, s} f_s + A^2_{t, s} f_s + f^\natural_{t, s}
$$
where
$$
f^\natural_{s, t} = \int_s^t (e^{\bf A}_{t, r} -\textrm{\emph{Id}}) B_r f_r \,dr + (e^{\bf A}_{t, s} -\textrm{\emph{Id}} - A^1_{t, s} - A^2_{t, s}) f_s.
$$
Using the bounds
$$
| e^{\bf A}_{t, r} -\textrm{\emph{Id}}| \lesssim_{\|{\bf A}\|} | t -  r |^{\gamma}, \quad \textrm{ and }\quad  | e^{\bf A}_{t, s}  -\textrm{\emph{Id}} - A^1_{t, s} - A^2_{t, s} | \lesssim_{\| {\bf A} \|} | t - s |^{3 \gamma},
$$
this allows to conclude that $|f^{\natural}_{t, s}| \lesssim | t - s |^{3\gamma}$, and that the path $f_\bullet$ is indeed the unique solution to equation \eqref{EqRoughLinearDrift}. \vspace{0.2cm}

   \item Bounded rough drivers have also been introduced previously in the work \cite{CoutinLejayLinear} of Coutin and Lejay, and studied in relation with the Magnus formula for what is called there the \emph{resolvent operator} $e^{\bf A}$. The above short proof of Theorem~\ref{ThmIntegrationLinearBounded} can be considered an alternative proof of the main result of section 3 in \cite{CoutinLejayLinear}. They also consider perturbed linear equations, with an a priori given drift of the more general form $C_{ts}$, rather than $\int_s^t B_rf_r\,dr$, with $C$ satisfying some regularity conditions. The pioneering work \cite{FdlPM} of Feyel-de la Pradelle-Mokobodzki is also closely related to these questions.
\end{enumerate}
\end{remark}

\medskip

\section{Unbounded rough drivers and rough linear equations}
\label{SectionUnboundedRoughDrivers}

The above results apply in the particular case where $\mcA$ is the Banach algebra of bounded operators on an Hilbert space $H$. We shall study in the remaining sections the integration problem 
\begin{equation}
\label{EqLinearEquationHilbert}
\delta f_{ts} = (A^1_{ts} + A^2_{ts}) f_s +  f^\natural_{ts},
\end{equation}
for a particular class of drivers $\bf A$ associated to a class of \textit{unbounded} operators on $H$, or other Banach spaces, with in mind the model case of the \textbf{rough transport equation}
$$
\delta f(\varphi)_{ts} = X_{ts}\,f_s(V^*\varphi) +  \bbX_{ts}\,f_s(V^*V^*\varphi) + f^\natural_{ts}(\varphi),
$$
where ${\bfX} = (X,\bbX)$ is an $\ell$-dimensional $\gamma$-H\"older rough path and $V = (V_1,\dots,V_\ell)$ is a collection of $\ell$ vector fields on $\RR^d$.

\subsection{Rough drivers}
\label{SubsectionRoughDrivers}

To make sense of this equation we need to complete the functional setting by the datum of a scale of Banach spaces $(E_n,|\cdot|_n)_{n\geq 0}$, with $E_{n+1}$ continuously embedded in $E_n$. For $n\geq 0$, we shall denote by $E_{-n}=E_n^*$ the dual space of $E_n$, equipped with its natural norm, 
$$
|e|_{-n} := \underset{\varphi\in E_n,\,|\varphi|_n\leq 1}{\sup}\; (\varphi,e), \quad\quad e\in E_{-n}.
$$
We require that the following continuous inclusions 
\begin{equation*}
E_n\subset \cdots\subset E_2\subset E_1\subset E_0 \\ 
\end{equation*}
hold for all $n\geq 2$. One can think of $n$ as quantifying the 'regularity' of elements of some test functions, with the elements of $E_n$ being more regular with $n$ increasing. Denote by $\|\cdot\|_{(b,a)}$ for the norm of a linear operator form $E_a$ to $E_b$. (Note that we use $(b,a)$ and not $(a,b)$ in the lower index for the norm.) We also assume the existence of a family $(J^\varepsilon)_{0<\varepsilon\leq 1}$ of operators from $E_0$ to itself such that the estimates
\begin{equation}
\label{EqApproximationProperties}
\|J^\varepsilon - \textrm{Id} \|_{(n,n+k)} \leq c\,\varepsilon^k, \quad \| J^\varepsilon \|_{(n+k,n )} \leq c\,\varepsilon^{- k}
\end{equation}
hold for all $n, k \geq 0$, for some positive constant $c$ independent of $\varepsilon$. For $\varphi\in E_0$, the elements $\varphi_\varepsilon := J^\varepsilon\varphi$ are in particular 'smooth', that is in the intersection of all the spaces $E_n$, for $n\geq 0$. 

Whenever we will work with Sobolev--like scales $E_n = W^{n,p}(\RR^d)$ ($p\ge 1$) we will take the operators $J^\varepsilon = (\textrm{I}-\varepsilon\triangle)^{-j_0}$, for $j_0$ big enough.

\smallskip

\begin{definition}
\label{DefnUnboundedRoughDriver}
Let $\frac{1}{3}<\gamma\leq \frac{1}{2}$ be given. An \textbf{\emph{unbounded $\gamma$-rough driver on the scales}} $(E_n,|\cdot|_n)_{n\geq 0}$, is a pair ${\bf A} = (A^1,A^2)$ of $2$-index maps, with 
\begin{equation}
\label{EqBoundednessAssumptionRoughDriver}
\begin{split}
&A^1_{ts}\in \textrm{\emph{L}}(E_n,E_{n-1}), \textrm{ for } n\in\{-0,-2\}, \\
&A^2_{ts}\in\textrm{\emph{L}}(E_n,E_{n-2}), \textrm{ for } n\in\{-0,-1\},
\end{split}
\end{equation}
for all $0\leq s\leq t<T$, which satisfies Chen's relations \eqref{EqChenRelations}, and such that $A^1$ is $\gamma$-H\"older and $A^2$ is $2\gamma$-H\"older. 
\end{definition}

\medskip

Equation \eqref{EqDefnSolutionRDE} below will make it clear that unbounded rough drivers can be thought of a some kind of (dual objects to some) multi-scales velocity fields, with two time scales given by $A^1$ and $A^2$. This is particularly clear on the following elementary example, where
$$
A^1_{ts} = \sum_{i=1}^\ell X^i_{ts}V_i, \quad \textrm{and } \quad A^2_{ts} = \sum_{j,k=1}^\ell \bbX^{jk}_{ts}V_jV_k,
$$
where ${\bfX} = (X,\bbX)$ is a weak geometric $\gamma$-H\"older rough path over $\RR^\ell$, with $\frac{1}{3}<\gamma\leq \frac{1}{2}$, and $V_1,\dots, V_\ell$ are regular enough vector fields on $\RR^d$ with values in $N\times N$ matrices, understood here as first order differential operators. On can also take $V_i^* = -V_i - \textrm{div}V_i$ in the above definition of a rough driver.

\smallskip

In a probabilistic setting, the above rough path $\bfX$ could be the (Stratonovich) Brownian rough path. More generally, one could take as first level $A^1_{ts}$ a semimartingale velocity field of Le Jan-Watanabe-Kunita type (or its dual), where noise (the above $X_{ts}$) and dynamics (the above $V_i$) cannot be separated from one another. It is shown in the forthcoming work \cite{BailleulRiedel} that these velocity fields can be lifted into a(n object very similar to a) rough driver under some mild conditions on the semimartingale structure. In the same spirit, Catellier has shown in \cite{Catellierregularisation, CatellierTransport} the interest for the study of rough transport equations with irregular drift, of considering velocity fields in $\RR^d$ given by the regularisation of some field $V$ (that may even be a distribution) along some irregular path $w$, such as a typical trajectory of a fractional Brownian motion with any Hurst index. One deals in this particular setting with the 'Young analogue' of the above $\gamma$-rough drivers, corresponding to $\frac{1}{2}<\gamma\leq 1$, with only one level, and given in this example by the formula
$$
A^1_{ts}(x) = \int_s^t V(x+w_u)\,du.
$$

\medskip

We denote by $\ep(e)$ the duality pairing between an element $\ep$ of $E_{-n}$ and an element $e$ of $E_n$, for any $n\in\ZZ$. The dual of a continuous operator $A$ from $E_{-a}$ to $E_{-b}$ is denoted by $A^*$; this is a continuous operator from $E_b$ to $E_a$. We can now make sense of equation \eqref{EqLinearEquationHilbert}. 

\begin{definition}
\label{DefnSolutionEquation}
An $E_{-0}$-valued \emph{\textbf{path}} $f_\bullet$ is said to \emph{\textbf{solve the linear rough differential equation}}
\begin{equation}
\label{EqRDE}
df_s = {\bf A}(ds)f_s
\end{equation}
if there exists an $E_{-2}$-valued $2$-index map $f^\natural$ such that one has
\begin{equation}
\label{EqDefnSolutionRDE}
\delta f_{ts}(\varphi) = f_s(\{A^{1,*}_{ts}+A^{2,*}_{ts}\}\varphi) + f^\natural_{ts}(\varphi)
\end{equation}
for all $0\leq s\leq t<T$ and all $\varphi\in E_2$, and the map $f^\natural_{ts}(\varphi)$ is $3\gamma$-H\"older, for each $\varphi\in E_2$.
\end{definition}

\smallskip

Let us insist on the fact that this notion of solution depends on the scale $(E_n)_{n\geq 0}$. We define, for all $0\leq s\leq t<T$, an  $E_{-1}$-valued $2$-index map $f^\sharp$ setting, for $\varphi\in E_1$,
\begin{equation*}
f^\sharp_{ts}(\varphi) := \delta f_{ts}(\varphi) - f_s(A_{ts}^{1,*} \varphi);
\end{equation*}
for $\varphi\in E_2$, one also has another expression for that quantity
$$
f^\sharp_{ts}(\varphi) = f_s(A^{2,*}_{ts}\varphi) + f^\natural_{ts}(\varphi).
$$

\bigskip

\subsection{A priori estimates}
\label{SubsectionAPrioriEstimates}

We show in this section that solutions $f_\bullet$ of equation \eqref{EqRDE} satisfy some a priori bounds that depend only on certains norms on the rough driver $A$ and on the uniform norm of $f$, when seen as a continuous path with values in some spaces of distributions. 

\smallskip

The next lemma shows that the map $f^\sharp$ is actually $2\gamma$-H\"older with values in a  space of distributions "of order 2" in a certain sense while $f^\sharp$ is only expected to be $\gamma$-H\"older from its definition. This \textit{gain of time regularity, traded against a loss of 'space regularity'} may well be one of our main contribution, despite its elementary nature. It leads to a similar (and even better) result for $f^\natural$, as expressed in Theorem~\ref{ThmGeneralRegularityGain}, that is the key technical result of this paper and which opens the road  to a quite complete theory of linear rough equation.

Recall the notations $\llparenthesis \cdot\rrparenthesis$ and $\llparenthesis \cdot\rrparenthesis_\gamma$ for the norms introduced in equations \eqref{EqDefnNormFunction} and \eqref{EqDefnNormMultiIndexMaps}. They depend implicitly on some positive parameter $\lambda$ that we shall tune as needed along the way. 
Given $\varphi\in E_0$, set 
\begin{equation*}
\begin{split}
K_1(\varphi) := \sup_{0 < \varepsilon \leq 1} \{& \llparenthesis\, f((A^1)^* J^\varepsilon\varphi)\,\rrparenthesis_\gamma
+ \lambda^{\gamma} \varepsilon \, \llparenthesis\, f( (A^2)^* J^\varepsilon \varphi) \,\rrparenthesis_{2\gamma}  \\ 
&+ \lambda^{2\gamma}\varepsilon^2\, \llparenthesis\, f^\natural (  J^\varepsilon\varphi) \,\rrparenthesis_{3\gamma} + \lambda^{-\gamma}\varepsilon^{-1} \,\llparenthesis\, f((\textrm{Id} - J^\varepsilon)\varphi) \,\rrparenthesis \}
\end{split}
\end{equation*}
and
\begin{equation*}
\begin{split}
K_2(\varphi) := \sup_{0 < \varepsilon \leq 1} \{&  \llparenthesis\, f((A^2)^* \varphi) \,\rrparenthesis_{2\gamma}  + \lambda^{\gamma}\varepsilon\,\llparenthesis\, f^\natural (  J^\varepsilon  \varphi\,) \,\rrparenthesis_{3\gamma} + \\ 
& \lambda^{-2\gamma}\varepsilon^{-2} \, \llparenthesis\, f((\textrm{Id} - J^\varepsilon) \varphi) \,\rrparenthesis + \lambda^{-\gamma} \varepsilon^{-1} \llparenthesis\, f((A^1)^*(\textrm{Id} - J^\varepsilon)\varphi) \,\rrparenthesis_\gamma  \}
\end{split}
\end{equation*}

\smallskip

\begin{lemma}
\label{LemmaGeneralImproveEstimates}
We have $\llparenthesis\, \delta f(\varphi) \,\rrparenthesis_\gamma \lesssim\; K_1(\varphi)$, and $\llparenthesis\,  f^\sharp(\varphi) \,\rrparenthesis_{2\gamma} \lesssim\; K_2(\varphi)$.
\end{lemma}

\medskip

\begin{proof}
We start by decomposing $\varphi\in E_1$ as $\varphi = \varphi_1 + \varphi_2$ with $\varphi_1 = J^\varepsilon \varphi$ and $\varphi_2 = (\textrm{Id}-J^\varepsilon) \varphi$ giving
$$
|\delta f_{ts}(\varphi)| \leq |\delta f_{ts}(\varphi_1)| + |\delta f_{ts}(\varphi_2)|.
$$
Take $t,s$ such that ${|t-s|\le \lambda}$. The second term is bounded above by 
$$
e^{-\lambda^{-1}t} \frac{|\delta f_{ts}(\varphi_2)|}{|t-s|^\gamma} \leq  2 \lambda^{-\gamma}\varepsilon^{-1}\llparenthesis\, \varphi_2 \rrparenthesis ,
$$ 
where  $\varepsilon = \lambda^{-\gamma} | t - s |^\gamma \in(0,1]$ if $|t-s|\le \lambda$.
The first term can be estimated using properties \eqref{EqApproximationProperties} and the defining identity \eqref{EqDefnSolutionRDE}:
$$
e^{-\lambda^{-1}t} \frac{|\delta f_{ts}(\varphi_1)|}{|t-s|^\gamma} \leq
\llparenthesis  f(A^{1,*} \varphi_1) \rrparenthesis_\gamma
+ \lambda^{\gamma} \varepsilon  \llparenthesis f( A^{2,*} \varphi_1) \rrparenthesis_{2\gamma}  
+ \lambda^{2\gamma}\varepsilon^2 \llparenthesis f^\natural (  \varphi_1) \rrparenthesis_{3\gamma}.
$$ 

\medskip

$\bullet$  This gives us an upper bound for $|\delta f_{ts}(\varphi)|$ depending on $\varepsilon$ and $\varphi_1, \varphi_2$ which we bound with
\begin{equation*}
\llparenthesis \delta f(\varphi) \rrparenthesis_\gamma \leq  K_1(\varphi) .
\end{equation*}

\smallskip

$\bullet$ One can proceed in a similar way to estimate $|f^\sharp_{ts}(\varphi)|$, starting from the inequality
\begin{equation*}
\begin{split}
| f^\sharp_{ts}(\varphi)| &\leq | \delta f_{ts}(\varphi_1)| + | f_s(A^{1,*}_{ts}\varphi_1)| + | f_s(A^{2,*}_{ts}\varphi_2)| + |f^\natural_{ts}(\varphi_2)|
\end{split}
\end{equation*}
from which get the upper bound on $|f^\sharp_{ts}(\varphi)|$ as a function of $K_2(\varphi)$.
\end{proof}

\medskip

The estimates proved in Lemma~\ref{LemmaGeneralImproveEstimates} are all we need to give an upper bound on the $3\gamma$-H\"older norm of $f^\natural$ in terms of $f$ itself. As this result will be our key tool for proving a number of results in different situations, we formulate it here in some generality. Given a Banach space B and a B-valued path $m_\bullet$, we set
$$
\llparenthesis\,m\,\rrparenthesis_{\gamma\,;\,\textrm{B}} := \sup_{\stackrel{0\leq s<t\leq T}{{|t-s|\le \lambda}}}e^{-\lambda^{-1}t}\,\frac{| \delta m_{ts} |_\textrm{B}}{|t-s|^\gamma},
$$
for $\gamma\geq 0$.
 Set 
$$
N_1({\bfA}) := \sup_{0<\varepsilon\leq 1} \{ \varepsilon\| J^\varepsilon A^{1,*}\|_{\gamma\,;\,\textrm{L}(\textrm{F},\textrm{F})} + \varepsilon^2\| J^\varepsilon A^{2,*}\|_{2\gamma\,;\,\textrm{L}(\textrm{F},\textrm{F})} \}
$$ 
and
\begin{equation*}
\begin{split}
N_2({\bfA}) := \sup_{0<\varepsilon\leq 1}  \{ &  \|A^{1,*} J^\varepsilon  A^{2,*}\|_{3\gamma\,;\,\textrm{L}(\textrm{F},\textrm{E})} + \varepsilon \|A^{2,*} J^\varepsilon  A^{2,*}\|_{4\gamma\,;\,\textrm{L}(\textrm{F},\textrm{E})} \\
&+ \varepsilon^{-1}  \|(\textrm{Id} - J^\varepsilon)A^{2,*}\|_{2\gamma\,;\,\textrm{L}(\textrm{F},\textrm{E})} + \|A^{2,*} A^{1,*}\|_{3\gamma\,;\,\textrm{L}(\textrm{F},\textrm{E})} \\
&+ \varepsilon^{-2} \|(\textrm{Id} - J^\varepsilon) A^{1,*}\|_{\gamma\,;\,\textrm{L}(\textrm{F},\textrm{E})} + \varepsilon^{-1}\|(A^{1,*}(\textrm{Id} - J^\varepsilon) A^{1,*}\|_{2\gamma\,;\,\textrm{L}(\textrm{F},\textrm{E})} \}
\end{split}
\end{equation*}

\smallskip

\begin{theorem}
\label{ThmGeneralRegularityGain} 
Let $\textrm{\emph{E}}\subset E_0$ and $\textrm{\emph{F}}\subset E_3$ be two Banach spaces of distributions such that  $N_1({\bfA})$ and $N_2({\bfA})$ are finite. Let  $\lambda \leq 1$ be sufficiently small so that $\lambda^{\gamma} N_1({\bfA}) \leq 1$. Then any solution of the rough linear equation \eqref{EqRDE} satisfies the a priori bound 
$$
\llparenthesis\, f^\natural\,\rrparenthesis_{3\gamma\,;\,F^*} \leq 4 \lambda^{-2\gamma} N_2({\bfA})\,\llparenthesis\,f\,\rrparenthesis_{E^*}. 
$$
\end{theorem}

\smallskip

\begin{proof}
Given $\varphi\in F$, note that
$$
\delta f_{tus} = f^\sharp_{us}(A^{1,*}_{tu}\varphi) + \delta f_{us}(A^{2,*}_{tu}\varphi),
$$
is $3\gamma$-H\"older by Lemma~\ref{LemmaGeneralImproveEstimates}, so the sewing Lemma~\ref{LemmaSewingMap} can be used, giving
\begin{equation*}
\begin{split}
\llparenthesis\, f^\natural(\varphi)\,\rrparenthesis_{3\gamma} 
&\lesssim \llparenthesis\,  f^\sharp (A^{1,*}\varphi) \,\rrparenthesis_{3\gamma} +  \llparenthesis\,  \delta f (A^{2,*}\varphi) \,\rrparenthesis_{3\gamma}\\
&\lesssim  \|  K_2(A^{1,*}\varphi) \|_\gamma +  \|  K_1(A^{2,*}\varphi)\|_{2\gamma}\\
\end{split}
\end{equation*}
Now
\begin{equation*}
\begin{split}
 \llparenthesis\, f^\natural\,\rrparenthesis_{3\gamma;\textrm{F}^*}
&\leq  2 \lambda^{-2\gamma}\llparenthesis\, f\,\rrparenthesis_{\textrm{E}^*} N_2({\bfA}) +  \lambda^{\gamma} \llparenthesis\, f^\natural\,\rrparenthesis_{3\gamma;\textrm{F}^*}  N_1({\bfA})
\end{split}
\end{equation*}
where we assumed  $\lambda \leq 1$, with no loss of generality, to get a simpler expression for the formula. 
Taking $\lambda\leq 1$ small enough so that
$
\lambda^{\gamma} N_1({\bfA}) \leq 1/2,
$
we obtain the  inequality
$$
\llparenthesis\, f^\natural \,\rrparenthesis_{3\gamma\,;\,\textrm{F}^*} \leq  4 \lambda^{-2\gamma} 
\llparenthesis\, f\,\rrparenthesis_{\textrm{E}^*} N_2({\bfA}).
$$
\end{proof}

\bigskip

Taking
$
\textrm{E} = E_0$, and $\textrm{F} = E_3$, 
and assuming that 
$$
A^{1,*}_{ts} \in \textrm{L}(E_1,E_0) \cap \textrm{L}(E_3,E_2)
,\quad
\text{and}\quad 
A^{2,*}_{ts} \in \textrm{L}(E_3 , E_1) \cap \textrm{L}(E_2,E_0)
$$ 
for all $0\leq s\leq t\leq T$, then 
$
N_1({\bfA})+N_2({\bfA}) \lesssim C_0^2
$ 
where
\begin{equation}
\label{EqDefnC0}
C_0 := 1+\|A^1\|_{\gamma\,;\,(-0,-1)}  + \|A^2\|_{2\gamma\,;\,(-0,-2)} +  \|A^1\|_{\gamma\,;\,(-2,-3)}  + \|A^2\|_{2\gamma\,;\,(-1,-3)}  < \infty .
\end{equation}

\smallskip

Note that it follows from Banach uniform boundedness principle that for a solution path $(f,f^\natural)$, we have
$$
\| f^\natural\|_{3\gamma\,;\,-3} \leq \| f^\natural\|_{3\gamma\,;\,-2} < \infty.
$$
Theorem~\ref{ThmGeneralRegularityGain} shows that one has an explicit upper bound for $\| f^\natural_{ts}\|_{3\gamma\,;\,-3}$ in terms of $ \llparenthesis f\rrparenthesis_{-0}$ and $C_0$ only. This fact is recorded in the following :

\begin{theorem}
\label{ThmRegularityGain} 
For any $\lambda \leq 1$ such that $\lambda^\gamma C_0^2 \lesssim 1$ we have the a priori bounds
\begin{equation}
\label{EqFundamentalAPrioriEstimate}
\max\{  \llparenthesis\, \delta f \,\rrparenthesis_{\gamma\,;\, - 1} , \llparenthesis\, f^\sharp \,\rrparenthesis_{2 \gamma\,;\, - 2}, \llparenthesis\, f^\natural \,\rrparenthesis_{3 \gamma\,;\, - 3} \} \lesssim_{\gamma} C_0^2\,\lambda^{-2\gamma}\,\llparenthesis\, f\,\rrparenthesis_{-0} .
\end{equation}
\end{theorem}

\medskip

\section{Integration of unbounded rough drivers in Hilbert spaces}
\label{SectionHilbertSpaceTheory}

We develop in this section the theory of integration of unbounded rough drivers in the Hilbert space $\LL^2(\RR^d)=E_0=E_{-0}$. We are able to give a rather complete theory a class of drivers that enjoys the same algebraic properties as the rough drivers ${\bfA} = (X\,V,\bbX\,V\,V)$ involved in the rough transport equation, when the vector fields $V=(V_1,\dots,V_\ell)$ are divergence free. These drivers are called conservative. A general existence result for the rough linear equation $df_s = {\bfA}(ds)f_s$, driven by conservative drivers is given in Section~\ref{SubsectionIntegrationConservative}. To prove uniqueness of solutions to such equations, we develop in section \eqref{SubsectionTensorization} a robust tensorization argument for a larger class of unbounded rough drivers that is the key to obtain some a priori bounds. These bounds imply uniqueness for rough linear equations driven by conservative drivers under a mild additional assumption, but they also lead to a complete $\LL^2$-theory of rough transport equations, as illustrated in section \eqref{SubsectionL2Transport}.

\smallskip

Note that working in a Hilbert space setting, we have the continuous inclusions
\begin{equation}
\label{EqDoubleInclusion}
E_n\subset \cdots \subset E_1\subset E_0 = E_{-0}\subset E_{-1}\subset \cdots\subset E_{-n}.
\end{equation}

\medskip

\subsection{Conservative drivers}
\label{SubsectionIntegrationConservative}

We start with the simple situation where the driver is \emph{conservative} according to the following definition.

\begin{definition}
A rough driver is said to be \textbf{\emph{conservative}} if we have 
\begin{equation*}
\begin{split}
&{\bf (i)}\quad A^{1,*}_{ts} +A^1_{ts} = 0, \quad \textrm{ on } E_1, \\
&{\bf (ii)}\quad A^{2,*}_{ts} + A^2_{ts} + A_{ts}^{1,*}  A^1_{ts} = 0, \quad \textrm{ on } E_2,
\end{split}
\end{equation*}
for all $0\leq s\leq t<T$.
\end{definition}
Notice that the conservative conditions make sense because of the above continuous inclusions \eqref{EqDoubleInclusion}. As an elementary example of conservative unbounded driver, take as Banach spaces $E_n$ the Sobolev spaces $W^{n,2}(\RR^d)$, with norm $|\varphi|_n := \sum_{k=0}^n |\nabla^k\varphi|_{\LL^2}$, and consider the unbounded rough driver $\bf A$  given by the formula
\begin{equation}
\label{EqRoughDriverRoughPath}
A^1_{ts} = X^i_{ts}V_i,  \quad\quad A^2_{ts} = \bbX^{jk}_{ts}V_jV_k,
\end{equation}
for some $\ell$-dimensional geometric $\gamma$-rough path, with $\frac{1}{3}<\gamma\leq 1$, and some \emph{divergence-free} vector fields $(V_i)_{i=1..\ell}$ on $\RR^d$, with the latter understood as first order differential operators. Condition {\bf (ii)} partly plays the role in our setting that the notion of weak geometric rough path plays in a rough paths setting. Also, the antisymmetry condition {\bf (i)} holds due to the fact that $V_i$ have null divergence, and condition {\bf (ii)} holds as a consequence of the weak geometric character of $\bfX$. Indeed, in this setting we have
\begin{equation*}
\begin{split}
&A^{1,*}_{ts} = X^i_{ts}(V_i)^* = - X^i_{ts} V_i = -A^1_{ts}, \\
&A^{2,*}_{ts} + A^2_{ts} =  \bbX^{jk}_{ts}(V_k^*V_j^*+V_j V_k) = \frac1 2  X^j_{ts}X^k_{ts} (V_kV_j +V_j V_k) = - A^{1,*}_{ts} A^1_{ts}
\end{split}
\end{equation*}
on $E_1$ and $E_2$ respectively. The boundedness assumptions \eqref{EqBoundednessAssumptionRoughDriver} that $A^1_{ts}$ and $A^2_{ts}$ need to satisfy hold if, for instance, the vector fields $V_i$ are $\mcC^2_b$. (Note that our setting is by no means restricted to working with vector fields. Working in the spaces $E_n = W^{kn,2}(\RR^d)$, we may take $V_i$ in the above formula for $\bfA$ the operators $V_i=W_i^k$, for some divergence-free vector fields $W_i$. Less trivial operators appear in a number of examples.)

\smallskip

A general existence result holds for equations driven by conservative rough drivers under very mild conditions on the functional setting.

\begin{theorem}
\label{ThmRDEHilbertGeneral}
Assume one can associate to the scale $(E_n)_{n\geq 0}$ a family of \emph{self-adjoint} smoothing operators $(J^\epsilon)_{0<\epsilon\leq 1}$, from $E_0$ to itself, satisfying the regularisation estimates \eqref{EqApproximationProperties}, and let $\bfA$ be a conservative unbounded $\gamma$-rough driver on the scales $(E_n)_{n\geq 0}$. Then for any $f_0\in\LL^2(\RR^d)$, there exists an $\LL^2(\RR^d)$-valued path $f_\bullet$, started from $f_0$, such that we have
\begin{equation*}
\delta f_{ts}(\varphi) = f_s(A^{1,*}_{ts}\varphi) + f_s(A^{2,*}_{ts}\varphi) + f^\natural_{ts}(\varphi)
\end{equation*}
for all $\varphi\in E_3$, with 
\begin{equation*}
| f_t|_0\leq | f_0|_0,
\end{equation*}
for all $ t\geq 0$, and we have, for each finite time horizon $T$,
\begin{equation}
\label{EqBoundRemainderX}
| f^\natural_{ts}(\varphi) | \lesssim_{C_0,{\bfA}, T,| f_0|_0}\,|\varphi|_3\,|t-s|^{3\gamma},
\end{equation}
for $0\leq s\leq t\leq T$.
\end{theorem}

\noindent The proof goes by approximating the unbounded rough driver $\bfA$ by bounded rough drivers  ${\bfA}^\epsilon$, and by using the theory developed in Section~\ref{SectionBoundedDrivers} to solve the equation
\begin{equation}
\label{EqApproximatedRDE}
\delta f^\epsilon_{ts}(\varphi) = f^\epsilon_s (A^{\epsilon,1,*}_{ts}\varphi) + f^\epsilon_s( A^{\epsilon,2,*}_{ts} \varphi) + f^{\epsilon,\natural}_{ts}(\varphi).
\end{equation}
The a priori bound \eqref{EqFundamentalAPrioriEstimate} is used to pass to the limit as  $\epsilon$ goes to $0$.

\medskip

\begin{proof}
More concretly, let 
$$
A^{\epsilon,1}_{ts} := J^\epsilon A^1_{ts}J^\epsilon,\quad\quad A^{\epsilon,2}_{ts} := J^\epsilon A^{2,a}_{ts}J^\epsilon - \frac1 2 A^{\epsilon,1}_{ts} A^{\epsilon,1,*}_{ts},
$$
with $J^\epsilon$ understood as a continuous operator from $E_0$ to $E_{k_0+2}$, for $k_0$ big enough, and $A^{2,a}_{ts} = \frac{1}{2}(A^2_{ts}-A^{2,*}_{ts})$ stands for the antisymmetric part of $A^{2}_{ts}$. Both $A^{\epsilon,1}$ and $A^{\epsilon,2}$ are bounded linear operators from $E_0$ to itself, and by construction ${\bfA}^\epsilon = (A^{\epsilon,1},A^{\epsilon,2})$ is a \emph{bounded} conservative $\gamma$-rough driver on the Banach algebra $\textrm{L}(E_0)$. So one can denote by $M^\epsilon_\bullet$ the solution path in $\textrm{L}(E_0)$ of the rough differential equation
$$
\delta M^\epsilon_{ts} = (A^{\epsilon,1}_{ts} + A^{\epsilon,2}_{ts})M^\epsilon_s + M^{\epsilon,\natural}_{ts},
$$
provided by the theory developed in Section~\ref{SectionBoundedDrivers}. Given $f_0\in E_0$, the path $f^\epsilon_\bullet = M^\epsilon_\bullet f_0$ solves the rough differential equation \eqref{EqApproximatedRDE}. Whereas Theorem~\ref{ThmIntegrationLinearBounded} provides an exponential control of the growth of the $E_0$-norm of $f^\epsilon_t$, with an exponent $\lambda^\epsilon$ that may go to $\infty$ as $\epsilon$ goes to $0$, the conservative character of ${\bfA}^\epsilon$ ensures the uniform bound
\begin{equation*}
| f^\epsilon_t|_0 = | f_0|_0.
\end{equation*}
Indeed, denoting by $(\cdot,\cdot)$ the scalar product on $E_0$, and setting $z^\epsilon_t = |f^\epsilon_t|^2_0$, one has
\begin{equation*}
\begin{split}
\delta z^\epsilon_{ts} &= 2(\delta f^\epsilon_{ts},f_s) + (\delta f^\epsilon_{ts}, \delta f^\epsilon_{ts}) \\
                                    &= 2((A^{\epsilon,1}_{ts}+A^{\epsilon,2}_{ts})f_s , f_s) + (f^{\epsilon,\natural}_{ts},f_s) \\
                                    &\quad\quad\quad\quad\quad\quad\quad\quad+ ((A^{\epsilon,1}_{ts}+A^{\epsilon,2}_{ts})f_s , (A^{\epsilon,1}_{ts}+A^{\epsilon,2}_{ts})f_s) + O(|t-s|^{3\gamma}) \\
                                    &= 2(A^{\epsilon,2}_{ts}f_s , f_s) + (A^{\epsilon,1}_{ts}f_s , A^{\epsilon,1}_{ts}f_s) + O(|t-s|^{3\gamma})
\end{split}
\end{equation*} 
by the conservative character of $A^{\epsilon,1}_{ts}$, which finally gives
$$
\delta z^\epsilon_{ts} = O(|t-s|^{3\gamma}),
$$
since the symmetric part of $A^{\epsilon,2}_{ts}$ is $-\frac{1}{2}A^{\epsilon,1,*}_{ts}A^{\epsilon,1}_{ts}$. The above equality shows that $z^\epsilon$ is constant.

\smallskip

Remark at that point that the constant 
\begin{equation*}
C^\epsilon_0 := \|A^{\epsilon,1}\|_{\gamma\,;\,(-0,-1)}  + \|A^{\epsilon,2}\|_{2\gamma\,;\,(-0,-2)} +  \|A^{\epsilon,1}\|_{\gamma\,;\,(-2,-3)}  + \|A^{\epsilon,2}\|_{2\gamma\,;\,(-1,-3)}  < \infty
\end{equation*}
is not only finite but also uniformly bounded above, independently of $\epsilon>0$. So we have, by Theorem~\ref{ThmRegularityGain}  an $\epsilon$-uniform upper bound on $\|f^{\epsilon,\natural}\|_{3\gamma\,;\,-3}$, of the form
\begin{equation}
\label{EqBondfSharpEpsilon}
\|f^{\epsilon,\natural}\|_{3\gamma\,;\,-3} \lesssim_{\gamma,T,| f_0|_0}\,1. 
\end{equation}
These bounds ensure in particular that for each $\varphi\in E_3$, the functions $f^\epsilon_\bullet(\varphi)$ form a bounded family of $\gamma$-H\"older real-valued paths, so it has a subsequence converging uniformly to some $\gamma$-H\"older real-valued function. Moreover, by weak-$\star$ compactness, the uniform bound \eqref{EqBondfSharpEpsilon} implies the existence of a sequence $(\epsilon_n)_{n\geq 0}$ converging to $0$, such that the sequence $f^{\epsilon_n}$ converges weakly-$\star$ in $\LL^\infty([0,T],E_0)$ to some limit $f\in\LL^\infty([0,T],E_0)$. In particular, for each $\varphi\in E_3$, the sequence $f^{\epsilon_n}(\varphi)$ converges weakly-$\star$ in $\LL^\infty([0,T],\RR)$ to $f(\varphi)$. As it has a uniformly converging subsequence, this shows the $\gamma$-H\"older character of each function $f(\varphi)$. 

Assuming $\varphi\in E_3$, it follows that one can pass to the limit in equation \eqref{EqApproximatedRDE} in the three terms involving $f^\epsilon_s$. The limit $f^\natural_{ts}(\varphi)$ is defined as a consequence, and the bound \eqref{EqBoundRemainderX} follows as a direct consequence of \eqref{EqBondfSharpEpsilon}.

\smallskip 

It is elementary to extend the above solution defined on $[0,T]$ to a globally defined solution satisfying the statement of the theorem.
\end{proof}

\bigskip

Rather than working with a general scale of spaces satisfying some ad hoc conditions, we shall set 
$$
E_n = W^{n,2}(\RR^d),
$$
for the remainder of this section on linear rough differential equation on Hilbert spaces. So we shall essentially be working from now on with rough drivers given by (at most) first order rough (pseudo-)differential operators.

\bigskip

\subsection{Tensorization}
\label{SubsectionTensorization}

In order to study the problem of uniqueness and further properties of solutions to general linear rough equations associated to unbounded rough drivers, we develop in this section a tensorization argument which can be seen as a rough version of the (differential) second quantisation functor in Hilbert spaces \cite{ReedSimon}, or the variables doubling method commonly used in the theory of transport equations and conservation laws after the pioneering work of Kruzkhov \cite{Kruzkhov}. As far as applications are concerned, we shall not restrict much our range in assuming that the rough drivers we are working with enjoy the following property. Given a bounded function $\phi\in W^{n_0,\infty}$, denote by $M_\phi$ the multiplication operator by $\phi$; it is a bounded operator from $E_0 = \LL^2(\RR^d)$ to itself.

\begin{definition}
\label{DefnSymmetricDriver}
An unbounded rough driver $\bfA$ is said to be \textbf{\emph{symmetric}} if the symmetric operators 
\begin{equation*}
\begin{split}
&{\bf (i)}\quad B^1_{{\bfA},ts}(\phi) = A^{1,*}_{ts}M_\phi  + M_\phi A^1_{ts} , \\
&{\bf (ii)}\quad B^2_{{\bfA},ts}(\phi) = M_\phi A^{2,*}_{ts}  +  A^2_{ts}M_\phi + A_{ts}^1 M_\phi A^{1,*}_{ts},
\end{split}
\end{equation*}
define quadratic forms 
$$
g\mapsto (g,B^i_{{\bfA},ts}(\phi)g)
$$
that are continuous on $E_0$, for all $0\leq s\leq t<T$, and for any test function $\phi\in W^{3,\infty}$.
\end{definition}

For rough drivers ${\bfA} = (A^1,A^2)$ of the form 
$$
A^1_{ts} = X^i_{ts} V_i, \quad\quad A^2_{ts} = {\bbX}^{jk}_{ts} V_jV_k,
$$
for some vector fields $V_i$ on $\RR^d$, and some weak geometric rough path $\bfX$, the operator $B^1(\phi)$ is the multiplication operator by $A^{1,*}_{ts}\phi$, and the second quadratic form is of the form
$$
(g,B^2_{{\bfA},ts}(\phi)g) = (g^2 , h^\phi_{ts}),
$$
for some explicit function $h^\phi_{ts}$ that turns $\bfA$ into a symmetric unbounded rough driver if the vector fields $V_i$ are $\mcC^2_b$ and $\phi\in W^{2,\infty}$. 

\bigskip 
 
Let now $\bfA$ be a symmetric unbounded rough driver in the scale $(E_n)_{n\geq 0}$, and let $f_\bullet$ be a solution to the linear rough differential equation 
$$ 
df_s = {\bfA}(ds)f_s.
$$
Consider $f^{\otimes 2} = f \otimes f$, defined on $\RR^d\times\RR^d$ by
$
f^{\otimes 2}(x,y) = f(x)f(y)
$; 
it satisfies the equation
$$ 
\delta f^{\otimes 2} = ( A^1  \otimes \mathbbm{I}+\mathbbm{I} \otimes A^1 ) f^{\otimes 2} + (A^2 \otimes \mathbbm{I}+\mathbbm{I} \otimes A^2 + \tau (A^1 \otimes A^1)) f^{\otimes 2} + f^{\otimes 2, \natural},
$$
where $\mathbbm{I}$ stands for the identity map and
$$
\tau (A^1 \otimes A^1)_{t s} = A^1_{t s} \otimes A^1_{t s}.
$$
Setting 
\begin{equation*}
\begin{split}
&\Gamma_{\bfA}^1 := A^1  \otimes \mathbbm{I}+\mathbbm{I} \otimes A^1, \\
&\Gamma_{\bfA}^2 := A^2 \otimes \mathbbm{I}+\mathbbm{I} \otimes A^2 + \tau (A^1 \otimes A^1),
\end{split}
\end{equation*}
and 
$$
\Gamma_{\bfA} := (\Gamma_{\bfA}^1, \Gamma_{\bfA}^2), 
$$ 
we have then 
$$
\delta \Gamma_{\bfA}^2 = (A^1 A^1) \otimes \mathbbm{I}+\mathbbm{I} \otimes (A^1 A^1) + A^1 \otimes A^1 + \sigma (A^1 \otimes A^1),
$$ 
where $\sigma (v \otimes w) = w \otimes v$ is the exchange of the two factors in the tensor so that 
$$
\sigma (A^1 \otimes A^1)_{t u s} = A^1_{us} \otimes A^1_{tu} + A^1_{tu} \otimes A^1_{us}.
$$ 
If follows that $\Gamma_{\bfA}$ satisfies Chen's relations \eqref{EqChenRelations}. Endowing $E_n^{\otimes 2}$ with its natural Hilbert space structure, we turn $\Gamma_{\bfA}$ into a $\gamma$-H\"older unbounded rough driver in the scale $(E_n^{\otimes 2})_{n \geq 0}$, for which $f^{\otimes 2}$ happens to be a solution, in this scale of spaces, of the linear rough equation
$$
df^{\otimes 2}_s = \Gamma_{\bfA}(ds)f^{\otimes 2}_s. 
$$

\smallskip 

The goal of the doubling variable method is to use the dynamics of $f^{\otimes 2}$ to gain some information of its behaviour near the diagonal of $\RR^d \times \RR^d$ in order to control the trace of $f^{\otimes 2}$ on the diagonal. 

Thinking to this goal let us define a scale of spaces $\mcE^\nabla_n$ of test functions on  $\RR^d \times \RR^d$. A function $\Phi \in W^{n,\infty}(\RR^d \times \RR^d)$ belongs to $\mcE^\nabla_n$ if $\Phi(x,y)=0$ when $|x-y|\ge 1$ and if the norm
$$
| \Phi |^\nabla_n :=\sup_{0\le k+l\le n}  \int_{\RR^d} \sup_{z\in\RR^d} |\partial_z^k \partial_w^l\Phi(z+w, z-w)|\, dw
$$ is finite. The support condition complicates the construction of suitable smoothing operators for such spaces. A detailed proof that this is possible can be found in~\cite{deya_priori_2016}.
Moreover we introduce a blow-up transformation $T_\varepsilon$ on test functions as follows:
$$
T_\varepsilon \Phi(x,y) =\varepsilon^{-d} \Phi(x_+ + \frac{x_-}{\varep}, x_+ - \frac{x_-}{\varep}) 
$$
where $x_\pm = \frac{x \pm y}{2}$ are coordinates parallel and transverse to the diagonal of $\RR^d \times \RR^d$. 
%
%
The adjoint of this transformation for the $\LL^2$ scalar product reads
$$
T^*_\varepsilon F(x,y) = F(x_+ + \varepsilon x_- , x_+ -\varepsilon x_-).
$$
We let  
$$
\Gamma^{\ast}_{\bfA,\varepsilon} := T_\varepsilon^{-1}\Gamma_{\bfA}^{\ast}T_\varepsilon.
$$ 
Recall the definition of the operators $ B^1_{{\bfA},ts}(\phi)$ and $ B^2_{{\bfA},ts}(\phi)$ used in the definition of a symmetric unbounded rough driver. We use below brackets $\langle\cdot,\cdot\rangle$ to denote the $\LL^2$ scalar product. Note that $f_s$ is an element of $E_{-0} = E_0 = \LL^2(\RR^d)$, so $f_s^2$ is in $\LL^1(\RR^d)$. 
 
\begin{definition}[Renomalizable drivers]
The unbounded $\gamma$-rough driver $\bfA$ is \textbf{\emph{renormalizable}} if $\{\Gamma_{\bfA,\varepsilon}\}_{\varepsilon>0}$  is a bounded family of unbounded $\gamma$-rough drivers in the scale of spaces $(\mcE^\nabla_n)_{n\geq 0}$.
\end{definition} 

The following lemma gives flesh to the expresion 'Renormalizable driver', such as defined here. It needs to be understood in the light of di Perna-Lions' work on transport equation \cite{diPernaLions} where a notion of renormalizable solution was first introduced. Definition~\ref{DefnrenormalisedSolutions} will make that parallel clear in our study of the $\LL^\infty$ theory for the rough transport equation. Under some additional assumption stated in Definition~\ref{DefnClosedDriver}, the present notion of renormalizable unbounded rough driver will provide in Theorem~\ref{ThmRDEHilbertGeneralClosed} a general uniqueness result.  

\begin{lemma}[Renormalisation]
\label{Lemmarenormalisation}
Let $\bfA$ be a symmetric rough driver and $f_\bullet$ be a solution of the equation 
$$
df_s = {\bfA}(ds)f_s,
$$ 
in the initial scale of spaces $(E_n)_{n\geq 0}$. If $\bfA$ is renormalizable 
then the $\LL^1(\RR^d)$-valued path $f^2_\bullet$ satisfies, for all $\phi \in W^{3,\infty}$, the equation
\begin{equation} 
\label{EqPrimitiveGronwall}
\delta f^2(\phi)_{ts} =  \langle f_s , (B^1_{{\bfA},ts}(\phi) + B^2_{{\bfA},ts}(\phi)) f_s\rangle  + O (| \phi |_{W^{3,\infty}} | t - s |^{3 \gamma}) .
\end{equation}
 \end{lemma}

\smallskip
 
\noindent By polarisation the product $fg$ satisfies an equation analogous to equation \eqref{EqPrimitiveGronwall} if both $f$ and $g$ are solutions of the equation $df_s = {\bfA}(ds)f_s$, in the scale $(E_n)_{n\geq 0}$.

\bigskip

\begin{proof}
Note that $f_\bullet^{\otimes 2}$ satisfies the equation
\begin{equation}
\label{EqFOtimes2}
\delta f^{\otimes 2}_\varepsilon (\Phi)_{ts} = f^{\otimes 2}_{\varepsilon,s}((\Gamma_{\bfA,\varepsilon}^{1,*})_{ts}\Phi) +f^{\otimes 2}_{\varepsilon,s} ( (\Gamma_{\bfA,\varepsilon}^{2,*})_{ts}\Phi) +  f^{\otimes 2,\natural}( T_\varepsilon\Phi)_{ts}
\end{equation}
for all smooth functions $\Phi$ where $f^{\otimes 2}_\varepsilon = T^*_\varepsilon f^{\otimes 2}$.
Note that if we show that $f^{\otimes 2}_\varepsilon$ is uniformly bounded in $\mcE^\nabla_{-0}$ then from the hypothesis that  $\{\Gamma_{\bfA,\varepsilon}\}_{\varepsilon>0}$ is a bounded family of unbounded $\gamma$-rough drivers in the scale of spaces $(\mcE^\nabla_n)_{n\geq 0}$ we also have
 $\|f^{\otimes 2,\natural}(T_\varepsilon \Phi)\|_{3\gamma}\lesssim | \Phi|^\nabla_3$. As a consequence the $3\gamma$-H\"older norm of the remainders $f^{\otimes 2,\,\natural}(T_\varepsilon \Phi)$ are bounded uniformly in  $\varepsilon$ for fixed $\Phi$. 

Equation \eqref{EqPrimitiveGronwall} will come from taking in equation \eqref{EqFOtimes2} some functions $\Phi$ of the form
$$
\Phi(x,y) = \psi(x-y)\,\phi\left(\frac{x+y}{2}\right),
$$
and by letting $\varepsilon$ tend to $0$, after checking that some $\varepsilon$-uniform estimates hold for the different terms in \eqref{EqFOtimes2}.

\smallskip
 Cauchy--Schwartz inequality provides the  bound
\begin{equation*}
\begin{split}
|f^{\otimes 2}_s (T_\varepsilon \Phi)| &= \left|\int_{\RR^d \times \RR^d} f_s(x_+ + \varepsilon x_-) f_s(x_+ - \varepsilon x_-) \Phi(x_+ +  x_-,x_+ -  x_-)\,dx_+ dx_- \right|  \\
&\leq \max_\pm \int_{\RR^d \times \RR^d} |f_s(x_+ \pm \varepsilon x_-)|^2  |\Phi(x_+ +  x_-,x_+ -  x_-)|\,dx_+ dx_- \\
&\leq  \max_\pm   \int_{\RR^d \times \RR^d} |f_s(x_+ \pm \varepsilon x_-)|^2  \sup_z |\Phi(z +  x_-,z -  x_-)|\,dx_+ dx_- \\
&\leq |f|_{\LL^2}^2 \int_{\RR^d } \sup_{z}| \Phi(z +  w,z -  w)|\,dw  \\
&\leq C |\Phi|_{\mcE^\nabla_0} |f_s|_{\LL^2}^2. 
\end{split}
\end{equation*}
which shows that $f^{\otimes 2}_\varepsilon$ is uniformly bounded in $\mcE^\nabla_{-0}$.
Now, given a positive constant $\delta$, the fact that for any smooth function $g$ which is $\delta$-close in $\LL^2(\RR^d)$ of $f_s$, we have
$$ 
|f^{\otimes 2}_s (T_\varepsilon \Phi)- g^{\otimes 2} (T_\varepsilon \Phi) | \lesssim 2\delta |f|_{\LL^2} + \delta^2,
$$ 
uniformly in $\varepsilon$, and
$$
\lim_{\varepsilon\to 0} g^{\otimes 2} (T_\varepsilon \Phi) = \int_{\RR^d }|g(x)|^2 \phi (x)\,dx,
$$
shows that 
$$
f^{\otimes 2}_s (T_\varepsilon \Phi) \underset{\varepsilon\to 0}{\longrightarrow} f_s^2(\phi). 
$$
We also have the convergence 
$$
g^{\otimes 2} \left((\Gamma_{\bfA}^1)^*_{ts} T_\varepsilon \Phi\right) = ( A^1_{ts} g \otimes g + g \otimes A^1_{ts} g, T_\varepsilon \Phi ) \underset{\varepsilon\to 0}{\longrightarrow} 2 \int_{\RR^d } g(x) (A^1_{ts}g)(x) \phi (x)\, dx,
$$
which we can rewrite as
$$
g^{\otimes 2} \left((\Gamma_{\bfA}^1)^*_{ts} T_\varepsilon \Phi\right) \underset{\varepsilon\to 0}{\longrightarrow} 2( g , A^{1,*}_{ts} M_\phi g ) = (g,B^1_{{\bfA},ts}(\phi)g).
$$
Using in addition the boundedness on $\LL^2(\RR^d)$ of the quadratic form associated to $B^1_{ts}(\phi)$, one can then send $g$ to $f_s$, in $\LL^2(\RR^d)$, in the above convergence result and conclude that
$$
f^{\otimes 2}_s ((\Gamma_{\bfA}^1)_{ts}^* T_\varepsilon \Phi) \to ( f_s , B^1_{{\bfA},ts}(\phi) f_s). 
$$
Similarly, the boundedness in $\mcE^\nabla_0$ of the family $((\Gamma_A^2)^*_{ts}  T_\varepsilon \Phi)_{0<\varepsilon\leq 1}$, together with the boundedness on $\LL^2(\RR^d)$ of the quadratic form associated with $B^2_{{\bfA},ts}(\phi)$, show that 
$$
f^{\otimes 2}_s ((\Gamma_{\bfA}^2)_{ts}^*\Phi^\epsilon) \to ( f_s , B^2_{{\bfA},ts}(\phi) f_s);
$$
equation \eqref{EqPrimitiveGronwall} follows, as we have the $\varepsilon$-uniform bound $|f^{\otimes 2,\,\natural}_{ts}(T_\varepsilon \Phi)| \leq  |\Phi|_{3}^\nabla \leq |\phi|_{W^{3,\infty}}|t-s|^{3\gamma}$.  
\end{proof} 

\bigskip

This result is sufficient to prove that rough linear equations driven by conservative drivers are unique if the driver is symmetric. 

\begin{corollary}
\label{CorClosedSymmCons}
Let $\bfA$ be a renormalizable symmetric conservative unbounded $\gamma$-rough driver in the scale of spaces $(E_n)_{n\geq 0}$. Then the rough linear equation
$$
df_s = {\bfA}(ds)f_s
$$
has a unique solution in $\LL^2(\RR^d)$, started from any initial condition $f_0\in E_0$; it satisfies $|f_t|_0 = |f_0|_0$, for all times $t$.
\end{corollary}

\smallskip

\begin{proof}
It suffices to notice that since the driver $\bfA$ is conservative, we have $B^1_{{\bfA},ts}({\bf 1}) = B^2_{{\bfA},ts}({\bf 1}) = 0$, so it follows from equation \eqref{EqPrimitiveGronwall} that any solution path $f_\bullet$ has constant $\LL^2$-norm, which proves the uniqueness claim. Existence was proved in Theorem~\ref{ThmRDEHilbertGeneral}.
\end{proof}

\bigskip

\subsection{A priori bounds for closed symmetric drivers}
\label{SymmetricDrivers}

One cannot use directly the renormalisation lemma to get some closed equation for $f^2$ when $\bfA$ is non-conservative. We need for that purpose to assume that the symmetric unbounded rough driver $\bfA$ enjoys the following property. Given that $\bfA$ is symmetric, recall the definition of its associated family of symmetric operators $B^1_{{\bfA},ts}(\phi)$ and $B^2_{{\bfA},ts}(\phi)$ on $E_0$, indexed by $(s,t)$ and $\phi\in W^{3,\infty}$, given in Definition~\ref{DefnSymmetricDriver}.

\begin{definition}
\label{DefnClosedDriver}
A \emph{\textbf{symmetric unbounded rough driver}} $\bfA$ in the scale of spaces $(W^{n,2})_{n\geq 0}$, is said to be \emph{\textbf{closed}} if there exists some unbounded rough driver ${\bf B}=(B^1,B^2)$ in the scale of spaces $(W^{n,\infty})_{n\geq 0}$, such that we have  
$$
\langle g,B^1_{{\bfA},ts}(\phi)g \rangle = \langle g^2,(B^1_{ts})^*\phi\rangle, \quad \textrm{  and  } \quad \langle g,B^2_{{\bfA},ts}(\phi)g\rangle = \langle g^2,(B^2_{ts})^*\phi\rangle,
$$
for all $g\in E_0$ and $\phi\in W^{3,\infty}$. 
\end{definition}

\smallskip

As an example, it is elementary to check that the unbounded rough driver ${\bfA} = (XV,\bbX VV)$ used in the rough transport equation
\begin{equation*}
\delta f_{ts} = X\,Vf_s +  \bbX\,V V f_s + f^\natural_{ts}
\end{equation*}
with some $\gamma$-H\"older weak geometric rough path ${\bfX} = (X,\bbX)$, is closed and symmetric if the vector fields $V = (V_1,\dots,V_\ell)$ are $C^3_b$, in which case $\bfB = \bfA$. Another interesting class of examples of closed symmetric unbounded rough driver in the scale of spaces $(E_n)_{n\geq 0}$, is provided by the lift to rough drivers of $\mcC^3_b$-semimartingale velocity fields, as given in the theory of stochastic flows. This kind of stochastic velocity fields appear in the study of Navier--Stokes equation. See the work \cite{BailleulRiedel} for a thorough study of stochastic flows from this point of view. 

\medskip

Building on this notion of closed driver, the following statement provides amongst other things an a priori estimate on solutions of rough linear equations that plays in this setting the role played in the classical setting by a priori estimates obtained by any kind of Gronwall--type argument. The crucial point here is that no such Gr\"onwall machinery was available so far in a rough path--like setting; despite its elementary nature, this result may well be one of our main contributions. 

\smallskip

\begin{theorem}
\label{ThmRDEHilbertGeneralClosed}
Let $\bfA$ be a renormalizable  \emph{closed symmetric} unbounded $\gamma$-rough driver on the scales $(W^{n,2}(\RR^d))_{n\in\NN}$. Let $\bfB$ be its associated driver, and assume that the inequality
\begin{equation}
\label{BoundPhi}
|B^{1,*}_{t0}{\bf 1} | \vee |B^{2,*}_{t0} {\bf 1} | \leq c_t
\end{equation}
holds for all times $t$, for some time-dependent mositive constant $c_t$ such that $e^{-\lambda t}c_t$ tends to $0$ as $t$ goes to infinity, for any positive parameter $\lambda$, so
$$
\llparenthesis\,c_\bullet\,\rrparenthesis < \infty.
$$
 Then, given any $f_0\in\LL^2(\RR^d)$, there is at most one $\LL^2(\RR^d)$-valued solution path $f_\bullet$ to the equation
\begin{equation}
\label{EqLinearRDEHilbert}
\delta f_{ts}(\varphi) = f_s (A^{1,*}_{ts}\varphi ) + f_s(A^{2,*}_{ts}\varphi) + f^\natural_{ts}(\varphi),
\end{equation}
and we have, for each finite time horizon $T$, 
\begin{equation}
\label{EqBoundRemainder}
| f^\natural_{ts}(\varphi) | \lesssim_{{\bfB},T, | f_0|_0}\,|\varphi|_3\,|t-s|^{3\gamma},
\end{equation}
for all $\varphi\in W^{3,2}$, and all $0\leq s\leq t\leq T$. It satisfies the upper bound
\begin{equation}
\label{EqAPrioriBoundNorm}
| f_t|_0\lesssim_{{\bfB}, t} | f_0|_0  .
\end{equation}
\end{theorem}

\smallskip
 
\begin{proof}
\noindent Let $f_\bullet$ be a solution to the rough linear equation \eqref{EqLinearRDEHilbert} in the scale of spaces $(W^{n,2}(\RR^d))_{n\in\NN}$. Since $\bfA$ is closed, the $\LL^1(\RR^d)$-valued path $f^2_\bullet$ happens to be a solution to the rough linear equation
\begin{equation*} 
\delta f^2(\phi)_{ts} = f^2_s((B^1_{ts})^*\phi) + f^2_s((B^2_{ts})^*\phi) + (f^2)^\natural_{ts}(\phi),
\end{equation*}
in the scale of spaces $(W^{n,\infty}(\RR^d))_{n\in\NN}$. Denote by $C_0^{\bf B}$ the finite constant associated to the unbounded rough driver $\bfB$, as defined by equation \eqref{EqDefnC0}, with $\bfB$ in the role of $\bfA$. It follows from the general a priori estimates on solutions of rough linear equations proved in Theorem~\ref{ThmGeneralRegularityGain}, and the fact that $f^2$ is in $\LL^1(\RR^d)$, that 
\begin{equation}
\label{EqBoundf2}
\begin{split}
\llparenthesis\, (f^2)^\natural \,\rrparenthesis_{3\gamma\,;\,(W^{3,\infty})^*} &\lesssim_{\gamma,\lambda}C_0^{\bf B}\,\llparenthesis\, f^2 \,\rrparenthesis_{(\LL^\infty)^*} 
\lesssim _{\gamma,\lambda} C_0^{\bf B}\,\llparenthesis\, f^2_\bullet \,\rrparenthesis_{\LL^1}.
\end{split}
\end{equation}
But since we have the identity 
$$
f^2_t({\bf 1}) = f_0^2({\bf 1} + (B^1_{ts})^*{\bf 1} + (B^2_{ts})^*{\bf 1}) + (f^2)^\natural_{t0}({\bf  1})
$$
and the bound \eqref{BoundPhi}, we also have the estimate
\begin{equation*}
\begin{split}
\llparenthesis\, f^2_\bullet \,\rrparenthesis_{\LL^1} = \llparenthesis\, f^2_\bullet({\bf 1}) \,\rrparenthesis &\lesssim_{\llparenthesis\,c_\bullet\,\rrparenthesis} \,| f_0 |_{\LL^2} + \llparenthesis\, (f^2)^\natural_{\bullet 0} \,\rrparenthesis_{(W^{3,\infty})^*}  \\
&\lesssim_{\llparenthesis\,c_\bullet\,\rrparenthesis} \, ( | f_0 |_{\LL^2} + \lambda^{3\gamma}\,\llparenthesis\, (f^2)^\natural \,\rrparenthesis_{3\gamma\,;\,(W^{3,\infty})^*} ).
\end{split}
\end{equation*}
(Note that $(f^2)^\natural$, in the right hand side of the above inequality, is seen as a 2-index function.) Together with the bound \eqref{EqBoundf2}, this gives the upper bound
$$
\llparenthesis\, f^2_\bullet \,\rrparenthesis_{\LL^1} \lesssim_{\llparenthesis\,c_\bullet\,\rrparenthesis,\,C_0^{\bf B}} \,| f_0 |_{\LL^2}
$$
for $\lambda$ small enough, which implies uniqueness.
\end{proof}

\medskip

\subsection{Rough transport equation}
\label{SubsectionL2Transport}

Building on Theorem~\ref{ThmRDEHilbertGeneralClosed}, one can give a complete $\LL^2$-theory of rough transport equations 
\begin{equation*}
\delta f_{ts} = X\,Vf_s +  \bbX\,V V f_s + f^\natural_{ts}
\end{equation*} 
driven by non-divergence-free vector fields $V_i$ of class $W^{3,\infty}$.

\medskip

\begin{lemma}
\label{LemmaGammaA}
Let $\bfX$ be a geometric $\gamma$-H\"older rough path on $\RR^\ell$, and $V_1,\dots,V_\ell$ be $W^{3,\infty}$ vector fields on $\RR^d$. Then the operator $\Gamma_{\bfA}$ associated with ${\bfA} = (XV, \bbX VV)$ is renormalizable in the scale of spaces $(\mcE^\nabla_n)_{n\geq 0}$.
\end{lemma}

\smallskip 

\begin{proof}
For a  geometric rough path $\bfX = (X,\bbX)$, the operator $\Gamma_{\bfA}$ takes the form $\Gamma_{\bfA} = (X  \Gamma_{V}^{1}, \bbX \Gamma_{V}^{2})$, with
\begin{equation*}
\begin{split}
&\Gamma_{V}^{1} := V  \otimes \mathbbm{I}+\mathbbm{I} \otimes V, \\
&\Gamma_{V}^{2} := V V \otimes \mathbbm{I}+\mathbbm{I} \otimes V V + 2 (V \otimes V) = \Gamma_{V}^{1} \Gamma_{V}^{1}.
\end{split}
\end{equation*}
So  it is enough to show that the adjoints of these operators satisfy, uniformly in $\varepsilon$, the inequalities
\begin{equation}
\label{eq:est-gamma-v}
|\Gamma_{V,\varepsilon}^{1,\ast} \Phi|^\nabla_{n} \lesssim |V|_{W^{n+1,\infty}} |\Phi|^\nabla_{n+1},
\qquad
|\Gamma_{V,\varepsilon}^{2,\ast} \Phi|^\nabla_{m} \lesssim |V|_{W^{m+2,\infty}}^2 |\Phi|^\nabla_{m+2}
\end{equation}
for $n=0,2,\,m=0,1$, for smooth test functions $\Phi$ where $\Gamma_{V,\varepsilon}^{j,\ast} = T_\varepsilon^{-1} \Gamma_{V}^{j,\ast}T_\varepsilon $ for $j=1,2$.

\smallskip

Write $V = v_k \partial_k$ where $(v_k)_{k=1,\dots,d}$ are the coefficients of the vector fields in the canonical basis $(\partial_k)_{k=1,\dots,d}$ of derivations; with these notations, we have
$$
V^\ast = - v_k \partial_k - d,
$$ 
where $d := \textrm{div} v$, is the divergence of the vector field $V$. For $\Gamma_{V}^{1,\ast} $ we have the representation
$$
\Gamma_{V}^{1,\ast} = -v^+_k \partial^+_k - v^-_k \partial^-_k - 2 d^+ ,
$$  
where we denote $\partial^\pm_k := (\partial_k \otimes 1) \pm (1 \otimes \partial_k)$ and, for a real-valued function $h$ on $\RR^d$,  by
$$
2 h^\pm(x,y) := h(x) \pm h(y)
$$ 
the symmetric and antisymmetric lifts of operators and functions from $\RR^d$ to $\RR^d \times \RR^d$. Moreover we let
$$
2 h^\pm_\varepsilon(x,y) := h(x_+ +\varepsilon x_-) \pm h(x_+ -\varepsilon x_-)
$$ 
the blowup of lifted functions according to the transformation $h^\pm_\varepsilon = T_\varepsilon^{-1} h^\pm T_\varepsilon$.
Note that
$$
\partial^+_k T_\varepsilon = T_\varepsilon \partial^+_k
\qquad
\partial^-_k T_\varepsilon = \varepsilon^{-1} T_\varepsilon \partial^-_k
$$
so that
$$
\Gamma_{V,\varepsilon}^{1,\ast} = T_\varepsilon^{-1} \Gamma_{V}^{1,\ast}T_\varepsilon = - v^+_{k,\varepsilon} \partial^+_k - \frac{ v^-_{k,\varepsilon}}{\varepsilon} \partial^-_k - 2 d^+_\varepsilon .
$$

The first estimate in \eqref{eq:est-gamma-v} follows from the inequalities
$$
| a^{+}_\varepsilon \nabla^{+} \Phi|^\nabla_{n}+| \varepsilon^{-1} a^{-}_\varepsilon \nabla^{-} \Phi|^\nabla_{n} \lesssim ( |a|_{W^{n,\infty}} + |\nabla a|_{W^{n,\infty}}) |\Phi|^\nabla_{n+1}.
$$ 
The second inequality in \eqref{eq:est-gamma-v} is obtained by noting that we have 
\begin{equation*}
\begin{split}
|\Gamma_{V,\varepsilon}^{2,\ast} \Phi|^\nabla_{m} =  |\Gamma_{V,\varepsilon}^{1,\ast}\Gamma_{V,\varepsilon}^{1,\ast} \Phi|^\nabla_{m} &\lesssim |V|_{W^{m+1,\infty}} |\Gamma_{V,\varepsilon}^{1,\ast} \Phi|^\nabla_{m+1} \\
&\lesssim |V|_{W^{m+1,\infty}} |V|_{W^{m+2,\infty}} | \Phi|^\nabla_{m+2}.
\end{split}
\end{equation*}

\end{proof}

\smallskip

\begin{theorem}
\label{ThmTransportL2}
Let $\bfX$ be a geometric $\gamma$-H\"older rough path on $\RR^\ell$, and $V_1,\dots,V_\ell$ be $W^{3,\infty}$ vector fields on $\RR^d$. Then the rough transport equation 
$$
\delta f_s = (X_{ts}\,V + \bbX_{ts}\,VV)f_s + f^\natural_{ts}
$$ 
is well-posed.
\end{theorem}

\smallskip

\begin{proof}
Notice first that the regularity assumption on the $V_i$ puts us in a position to use the a priori bounds for symmetric closed drivers stated in Theorem~\ref{ThmRDEHilbertGeneral}, with $\bfA$ in the role of $\bfB$. So uniqueness is a direct consequence of the a priori bound \eqref{EqAPrioriBoundNorm} in Theorem~\ref{ThmRDEHilbertGeneralClosed}.

\smallskip

Let now $f_0\in E_0$ be given. We prove the existence of a solution path to rough transport equation started from $f_0$, by a classical approximation-compactness argument, relying in a crucial way on the a priori bound \eqref{EqAPrioriBoundNorm} on the $\LL^2$-norm of the solution to the approximate problem, and on the uniform estimate \eqref{EqBoundRemainderX} for the remainder. 

\smallskip

Fix a finite time horizon $T$.
Given that $\bfX$ is geometric, let $({\bfX}^\epsilon)_{0<\epsilon\leq 1}$ be a family of rough path lifts of smooth paths which converge to $\bfX$ is a rough paths sense over the time interval $[0,T]$. Let also $(V^\epsilon)_{0<\epsilon\leq 1}$ be a family of smooth vector fields that converge to $V$ in $W^{3,\infty}$, and let $({\bf A}^\epsilon)_{0<\epsilon\leq 1}$ be their associated rough driver, defined by formula \eqref{EqRoughDriverRoughPath} with ${\bfX}^\epsilon$ and $V^\epsilon$ in place of $\bfX$ and $V$ respectively. One can choose $({\bfX}^\epsilon)_{0<\epsilon\leq 1}$ in such a way that the constant $C_0^\epsilon$ associated with ${\bf A}^\epsilon$ by formula \eqref{EqDefnC0} satisfies the inequality $C_0^\epsilon\lesssim C_0$, independently of $0<\epsilon\leq 1$. Given the smooth character of the vector fields $V^\epsilon$, one can solve uniquely the transport equation
\begin{equation*}
\delta f^\epsilon_{ts}(\varphi) = f^\epsilon_s(V^{\epsilon,\ast} \varphi)\,X^\epsilon_{ts} + f^\epsilon_s(V^{\epsilon,\ast}V^{\epsilon,\ast}\varphi)\,\bbX^\epsilon_{ts} + f^{\epsilon, \natural}_{ts}(\varphi), \quad\quad \textrm{for } \varphi\in E_2,
\end{equation*}
by the method of characteristics, as the above equation is actually equivalent to the partial differential equation
$$ 
\frac{df^\epsilon_t(\varphi)}{dt} = f^\epsilon_t(V^{\epsilon,\ast}\varphi)\,\frac{dX^\epsilon_t}{dt}.
$$
The solutions of this problem satisfy the uniform estimates
\begin{equation*}
| f^\epsilon_t |_0 \lesssim_{C_0,T}\,| f_0 |_0,
\end{equation*}
for all $0\leq t\leq T$, as a consequence of \eqref{EqAPrioriBoundNorm}, and we also have the uniform bound
\begin{equation*}
\sup_{0<\epsilon\leq 1}\;\| f^{\epsilon, \natural} \|_{3 \gamma\,;\, - 3} \lesssim_{C_0, \gamma, T,|f_0|_0} 1,
\end{equation*}
by \eqref{EqBoundRemainderX}. These two a priori estimates are all we need to finish the proof of the theorem following word by word the end of the proof of Theorem~\ref{ThmRDEHilbertGeneralClosed}. 
\end{proof}

\medskip

It is perfectly possible to extend the present theory to deal with rough linear equations with a {\bf drift}
$$
df_s = Wf_s ds + {\bf A}(ds)f_s,
$$
where $W\in\textrm{L}(E_{-0},E_{-2})$, such as the Laplacian operator in the $W^{n,2}(\RR^d)$ scale of spaces. We refrain from giving the details here as this is not our main point and this does not require the introduction of new tools or ideas. This provides an alternative road to some of the results of \cite{DiehlFrizStannat} in a slightly different setting.

\bigskip

\section{The $\LL^\infty$ theory of rough transport equations}
\label{SectionLInftyTheory}

We develop in this section an $\LL^\infty$ theory of the rough transport equation 
\begin{equation}
\label{EqTransportEq}
\delta f_{ts} = X\,Vf_s +  \bbX\,V V f_s + f^\natural_{ts}
\end{equation}  
and prove its well-posed character under the assumption that the vector fields be $\mcC^3$. We show for that purpose that all solutions are renormalised solutions, in the sense of Di~Perna--Lions, which classically leads to uniqueness and stability results in that setting. 

\medskip

\subsection{A priori estimates and existence result}
\label{SubsectionAPrioriEstimatesLInfty}

For developing that $\LL^\infty$ theory, we shall be working in the scale of Sobolev spaces 
$$
E_n = W^{n,1}(\RR^d),\quad\textrm{for } n\geq 0,
$$
with norm denoted by $|\cdot|_n$, and in which one has regularising operators $(J^\ep)_{0<\epsilon\leq 1}$ for which estimates \eqref{EqApproximationProperties} hold. Our minimal regularity assumptions on the vector fields will be the existence of a positive constant $C_1$ such that the inequalities
\begin{equation}
\label{EqContinuityConditionsVStar}
|V_i^*\varphi|_0 \leq C_1 |\varphi|_1, \quad\quad |V_i^*V_j^*\varphi|_0 \leq C_1 |\varphi|_2
\end{equation}
hold for all $1\leq i,j\leq\ell$. These conditions hold for instance if the vector fields $V_i$ and $(V_iV_j)$ are all $\mcC^1_b$; we write here $(V_iV_j)$ for $(DV_j)(V_i)$. One proves the following existence result by proceeding exactly as in the proof of Theorem~\ref{ThmTransportL2}, using the a priori $\LL^\infty$-estimate
$$
| f_t^\ep |_{\LL^\infty} = |f_0^\ep|_{\LL^\infty},
$$
for the regularised equation, and using Theorem~\ref{ThmRegularityGain} to get an $\epsilon$-uniform control on $|f^{\epsilon,\natural}|_{3\gamma\,;\,-3}$, in terms of $\bfX$ and $|f_0^\ep|_{\LL^\infty}$ only. It holds in particular if $V$ is $\mcC^2_b$.

\smallskip 
 
\begin{theorem}[Existence for rough transport equations -- $\LL^\infty$ setting]
\label{ThmExistenceTransportEquationLInfty}
Under the continuity assumptions \eqref{EqContinuityConditionsVStar} on the vector fields $V_i$, for any $f_0\in\LL^\infty(\RR^d)$, there exists an $\LL^\infty(\RR^d)$-valued path $(f_t)_{t\geq 0}$, started from $f_0$, satisfying the equation
$$
\delta f_{ts}(\varphi) = f_s(V^*\varphi)\,X_{ts} + f_s(V^*V^*\varphi)\,\bbX_{ts} + f^\natural_{ts}(\varphi)
$$
for all $\varphi\in E_3$, and the bound
$$ 
\sup_{t\geq 0} | f_t|_{\LL^\infty(\RR^d)}\leq | f_0|_{\LL^\infty(\RR^d)},
$$ 
with a remainder $f^\natural(\varphi)$ controlled by 
\begin{equation}
\label{EqBoundRemainder-bis}
|f^\natural_{ts}(\varphi) | \lesssim_{C_1,{\bfX}, T, |f_0|_{\LL^\infty}} \, |\varphi|_3\,|t-s|^{3\gamma},
\end{equation}
for $0\leq s\leq t\leq T$. 
\end{theorem}

\bigskip

\subsection{Renormalised solutions, uniqueness and stability}
\label{SubsectionrenormalisedSolutions}

To proceed one step further, we show that a mild strengthening of the regularity conditions imposed on the vector fields $V_i$ suffices to guarantee that all bounded solutions to the transport equation \eqref{EqTransportEq} are actually renormalised solution, in the sense of the following definition. 

\begin{definition}
\label{DefnrenormalisedSolutions}
A solution $f_\bullet$ to the transport equation \eqref{EqTransportEq} in the scales $(E_n)_{n\geq 0}$ is said to be a \emph{\bf renormalised solution} if for any function $H : \RR\rightarrow\RR$, of class $\mcC^3_b$, the path $h_\bullet = H\circ f_\bullet$ is also a solution to equation \eqref{EqTransportEq} in the same scale $(E_n)_{n\geq 0}$.
\end{definition}

As expected, this property will lead below to uniqueness and stability results.

\begin{theorem}
\label{ThmRenomalizedSolutions}
Assume the vector fields $V_i$ are $\mcC^3_b$. Then every solution to the transport equation \eqref{EqTransportEq}, bounded in $ \LL^\infty(\RR^d)$, is a renormalised solution.
\end{theorem}
 
\medskip
 
\begin{proof}  
The renormalisation Lemma~\ref{Lemmarenormalisation} can be stated in the $\LL^\infty$ setting by chosing a slightly different scale $(\mcF^\nabla_n)_n$ of spaces of test functions. We say that $\varphi \in \mcF^\nabla_n \subseteq L^1(\RR^d \times \RR^d)$ if the norm
\begin{equation}
\label{EqNormNablaLInfty} 
| \varphi |^\nabla_n := \sup_{0 \leqslant k + \ell \leqslant n} \int \int |(\nabla^+)^k (\nabla^-)^\ell \varphi (x, y) |\,dx dy  
\end{equation}
is finite and the support of $\varphi$ is contained in the set $\{(x,y) : |x-y|\le 1\}$. Compare this $\LL^1$ space  with the $\LL^\infty$ space used in Section~\ref{SubsectionTensorization}. Identity \eqref{EqPrimitiveGronwall} holds in that case with for functions $\phi\in W^{3,1}(\RR^d)$, with an $O(\cdot)$ term involving the $W^{3,1}$-norm of $\phi$ rather than its $W^{3,\infty}$-norm, as the proof of Lemma~\ref{Lemmarenormalisation} works verbatim, provided we can prove that $\Gamma_{\bfA}$ is an unbounded rough driver in the scale of spaces $(\mcF_n^\nabla)_{n\geq 0}$ associated with the norm \eqref{EqNormNablaLInfty}.  (Note that we have in that case $| f_s^{\otimes 2}(T_\varepsilon \Phi)| \leq |\Phi|_0^\nabla | f_s|_{\LL^\infty}^2 = |\Phi|_0^\nabla | f_s|_{-0}^2$.)

\medskip

The proof that $\Gamma_{\bfA}$ is a unbounded rough driver renormalizable in the scale $(\mcF^\nabla_n)_n$ follows the same pattern as the proof given in Section~\ref{SubsectionL2Transport}. We invite the reader to complete the details.

\smallskip

So it follows from the renormalisation lemma that if $f, g$ are two solutions the above argument also goes through and shows that $fg$ is also a solution, so any power $f^n$ of $f$ is also a solution, with a size of the remainder that depends only on the $\LL^{\infty}$ norm of $f^n$. By linearity the result can be extended to any polynomial of $f$, and by density to any continuous function $H (f)$, with $H$ defined on the interval $[- \| f \|_{\infty}, \| f \|_{\infty}]$. 
\end{proof}
 
\medskip

We can actually improve slightly this condition and require only a weak integrability for the third derivative of $V$; it provides a significant strengthening of the previous statement when the vector fields $V_i$ are divergence-free, giving some analogue of the traditional di~Perna--Lions conditions in the classical setting. Note that we do not have uniqueness of solutions for the associated rough differential equation under the conditions below.

\begin{theorem}
\label{ThmRenomalizedSolutionsWeak}
Assume that $V \in \mcC^2_b$, $\nabla^3 V \in L^1$ and $\mathrm{div} V \in \mcC^2_b$. Then every solution to the transport equation \eqref{EqTransportEq}, bounded in $ \LL^\infty(\RR^d)$, is a renormalised solution.
\end{theorem}
 
\medskip
 
\begin{proof}  
In the proof of the renormalisation Lemma~\ref{Lemmarenormalisation} we can use directly the general a priori estimate stated in Theorem~\ref{ThmGeneralRegularityGain} applied to $\Gamma_{\bfA}$ with  $F=\widetilde \mcE^\nabla_3$ and $E= \mcF^\nabla_0$ -- note the choice of function space for $F$. Here $\widetilde \mcE^\nabla_n$ are spaces of test functions,  with norms modelled on $\LL^\infty$ like  $\mcE^\nabla_n$ and support contained in the set $\{(x,y): |x-y|\le 1\}$ but with a small change given by an additional averaging over the auxiliary variable $\tau$ and a weight:
\begin{equation*}
| \varphi |^\nabla_n := \sup_{0 \leqslant k + \ell \leqslant n} \int  \sup_{x\in\RR^d} [\int_0^1 d\tau |(\nabla^+)^k (\nabla^-)^\ell \varphi (x-\tau w, x+(1-\tau)w) |]\, (1+|w|) dw 
\end{equation*}
the reason of which will be clear below.
In this case we can  show that
$$
N_1(\Gamma_{\bfA,\varepsilon}) \lesssim (1+ |V|_{\mcC^2_b})^2
$$
while
$$
N_2(\Gamma_{\bfA,\varepsilon}) \lesssim (1+|V|_{\mcC^2_b}+|\nabla^3 V|_{\LL^1}+|\mathrm{div} V|_{\mcC^2_b})^3.
$$
Indeed apart from many contributions which can be estimated as in the $\LL^2$ or as in the previous theorem, a difficult term come form the estimation of norms like $|\Gamma_{V,\varepsilon}^\ast \Gamma_{V,\varepsilon}^\ast  \Gamma_{V,\varepsilon}^\ast  |_{\textrm{L}(\textrm{F},\textrm{E})}$ of which the most singular contribution is given by $|\Gamma_{V,\varepsilon} \Gamma_{V,\varepsilon}  \Gamma_{V,\varepsilon}  |_{\textrm{L}(\textrm{F},\textrm{E})}$. In this norm the contribution that requires more regularity to $V$ is due to  the first two vector fields $\Gamma_{V,\varepsilon}$ acting simultaneously on the third one giving terms of the form $\varepsilon^{-1}|v_\varepsilon^+ v_\varepsilon^+ (\nabla^2 v)^-_\varepsilon \nabla^- |_{\textrm{L}(\textrm{F},\textrm{E})}$ and easier ones. Now expanding $(\nabla^2 v_\varepsilon)^-(x,y)= \varepsilon \int_0^1 d\tau (\nabla^3v)(x+\varepsilon \tau(y-x)) (y-x)$ we get
$$
\varepsilon^{-1}|v^+_\varepsilon v^+_\varepsilon (\nabla^2 v)^-_\varepsilon \nabla^- \Psi|_E =\varepsilon^{-1}  \int \int |(v^+_\varepsilon v^+_\varepsilon (\nabla^2 v)^-_\varepsilon  \Psi) (x, y) | dx dy  $$
$$
 \lesssim  |V|_{\LL^\infty}^2 \int_0^1 d\tau \int \int |\nabla^3 v(x+\varepsilon \tau w)| |\Psi(x,x+w)| |w| dx dw 
 $$
$$
 \lesssim  |V|_{\LL^\infty}^2 \int_0^1 d\tau \int dx |\nabla^3 v(x)| \int dw  |\Psi(x-\tau w,x-\tau w+w)| |w|
 $$
$$
 \lesssim  |V|_{\LL^\infty}^2 |\nabla^3 V|_{\LL^1} \sup_{x\in\RR^d} \int_0^1 d\tau \int dw  |\Psi(x-\tau w,x-\tau w+w)| (1+|w|)
 $$
$$
 \lesssim  |V|_{\LL^\infty}^2 |\nabla^3 V|_{\LL^1}  |\Phi|_{F}.
 $$
Granted the bounds on $N_1(\Gamma_{\bfA})$ and $N_2(\Gamma_{\bfA})$ the proof continues as the proof of the previous theorem and gives the renormalisation result.
\end{proof}
 
\medskip

As expected, Theorem~\ref{ThmRenomalizedSolutions} on renormalised solutions to the transport equation \eqref{EqTransportEq} comes with a number of important consequences, amongst which is an equivalent of the missing Gronwall lemma, as given by the a priori estimate \eqref{EqGronwall} below.

\begin{theorem}
\label{ThmUniqueness} 
Assume the vector fields $V_i$ are $\mcC^3_b$. \vspace{0.1cm}
\begin{enumerate}
   \item {\bf Uniqueness --} Given an initial condition in $\LL^\infty(\RR^d)$, there exists a unique solution to the transport equation which remains bounded in $\LL^\infty(\RR^d)$.  \vspace{0.2cm}
   \item {\bf Stability --} Let the time horizon $T$ be finite. Let $(V^{(n)}_i)_{n\geq 0}$, $i=1..\ell$ and $(f^{(n)}_0)_{n\geq 0}$ be a sequence of approximating sequences with $V^{(n)}_i$ converging to $V_i$ in $\mcC^3_b$, and $f^{(n)}_0$ converging to $f_0$ in $\LL(\RR^d)$. Let also $({\bfX}^{(n)})_{n\geq 0}$ be a sequence of weak geometric $\gamma$-rough paths above smooth paths, that converge in a rough paths sense to $\bfX$, over the time interval $[0,T]$. Then the solution paths $f^{(n)}_\bullet$ to the transport equation associated with ${\bfX}^{(n)}, V^{(n)}$ and $f^{(n)}_0$, converge weakly-$\star$ in $\LL^\infty([0,T],\LL^\infty(\RR^d))$, and in $\LL^1_\textrm{\emph{loc}}([0,T],\LL^\infty(\RR^d))$, to $f_\bullet$.
\end{enumerate}
\end{theorem}

\medskip

\begin{proof}
{\bf Uniqueness --} We follow the same pattern of proof as that of Theorem~\ref{ThmRDEHilbertGeneralClosed}. Let $f_\bullet$ and $f'_\bullet$ be two solution paths to equation \eqref{EqTransportEq}, bounded in $\LL^\infty(\RR^d)$, and started from the same initial condition. Let $H : \RR\rightarrow\RR$ be a non-negative function, of class $\mcC^3_b$, null at $0$ and positive elsewhere. Define the path 
$$
h_\bullet = H(f_\bullet - f'_\bullet);
$$
it is also a positive solution to the transport equation \eqref{EqTransportEq} under the above regularity assumptions on the vector fields $V_i$, since all solutions are renormalised solution, by Theorem~\ref{ThmRenomalizedSolutions}. Set $\psi(x) = (1+|x|^2)^{-k_0}$, for $x\in\RR^d$, and some exponent $k_0>d$. That function satisfies 
\begin{equation*}
\begin{split}
&\textrm{div}(\psi V) = -V\psi - (\textrm{div}V)\psi, \\
&\textrm{div}(\textrm{div}(\psi V)\,V) = V^2\psi + (\textrm{div}V)V\psi + ((V\textrm{div}V)+(\textrm{div}V)^2)\psi  
\end{split}
\end{equation*}
with 
\begin{equation*}
|V\psi| \,\vee\, |V^2\psi + (\textrm{div}V)V\psi| \lesssim \psi
\end{equation*}
as a consequence of the $\mcC^1_b$ character of the vector fields $V_i$. 
Define the scale of spaces 
$$
E_n^\psi := \{\varphi = \psi \phi\,;\,\phi \in \LL^\infty\} 
$$
with norm 
$$
|\varphi|_{E_n^\psi} := |\phi|_{W^{n,\infty}}
$$ 
It is not difficult to check that $(V X, VV \mathbb{X})$ is a $\gamma$-rough driver also in this scale of spaces. In this case however we have
$$
|h_t(\varphi)| = |h_t ( \psi \phi)| \leq |\phi|_{\LL^\infty} |h_t( \psi)| = |\varphi|_{E^\psi_0}|z_t|
$$
where 
$$
z_t := h_t(\psi), 
$$
and so
$$
\llparenthesis h_\bullet\rrparenthesis_{(E_0^\psi)^*}\leq \llparenthesis z_\bullet \rrparenthesis
$$
By the general a priori estimates we have that there exists $\lambda$ and constants depending on $A$ such that
$$
\llparenthesis h^\natural \rrparenthesis_{3\gamma\,;\,(E_3^\psi)^*} \lesssim \llparenthesis z_\bullet \rrparenthesis
$$
But now
$$
z_t = z_0 +  h_{0}(V^* \psi)X_{0,t}+  h_{0}(V^*V^* \psi)\mathbb{X}_{0,t} + h^\natural_{0,t}(\psi)
$$
so
$$
\llparenthesis z_\bullet \rrparenthesis \leq |z_0| (1+ \llparenthesis X_{0, \bullet} \rrparenthesis+ \llparenthesis \mathbb{X}_{0,\bullet} \rrparenthesis) + \llparenthesis h^\natural_{0,\bullet} \rrparenthesis_{(E_3^\psi)^*}
$$
and since $  \llparenthesis h^\natural_{\bullet 0} \rrparenthesis_{(E_3^\psi)^*} \leq \lambda^{3\gamma} \llparenthesis h^\natural \rrparenthesis_{3\gamma\,;\,(E_3^\psi)^*} \leq \lambda^{3\gamma}  \llparenthesis z_\bullet \rrparenthesis$, with $h^\natural$ considered as a 2-index map in its second occurence, we have for $\lambda$ small enough
\begin{equation}
\label{EqGronwall}
\llparenthesis z_\bullet \rrparenthesis \leq 2 |z_0| (1+ \llparenthesis X_{0, \bullet} \rrparenthesis+ \llparenthesis \mathbb{X}_{0, \bullet} \rrparenthesis);
\end{equation}
so $z_t = 0$, for all $t\ge 0$, if $z_0 = 0$.

\bigskip
  
{\bf Stability --} Denote by ${\bfX}^{(n)}$ a smooth rough path converging to $\bfX$ in the rough paths metric, and by $V^{(n)}_i$ a sequence of vector fields converging to $V_i$ in $\mcC^3_b$. Let $f^{(n)}_0$ be a smooth function converging to $f_0$ in $(\LL^1)^*$. One solves the transport equation associated with ${\bfX}^{(n)}$ and $f^{(n)}_0$, using the elementary method of characteristics as the vector fields $V^{(n)}_i$ are sufficiently regular. It is elementary to use the uniform bound
$$
\left|f^{(n)}_t\right|_{\LL^\infty(\RR^d)} \leq \left|f^{(n)}_0\right|_{\LL^\infty(\RR^d)} \leq C<\infty,
$$ 
and the uniform a priori bound on $|f^{(n),\natural}|_{3\gamma\,;\,-3}$ provided by Theorem~\ref{ThmRegularityGain} and the convergence of ${\bfX}^{(n)}$ to $\bfX$, and $f^{(n)}_0$ to $f_0$, to get the existence of a subsequence of $f^{(n)}_\bullet$ converging weakly-$\star$ in $\LL^\infty([0,T];\LL^\infty(\RR^d))$ to some solution of the transport equation \eqref{EqTransportEq}, bounded in $\LL^\infty(\RR^d)$. Since this solution is unique, as proved above, the whole sequence $f^{(n)}_\bullet$ converges weakly-$\star$ to $f_\bullet$ in $\LL^\infty([0,T];\LL^\infty(\RR^d))$. As the same conclusion holds for $(f^{(n)}_\bullet)^2$ and $f^2$, by the renormalisation property, we classically get the convergence in $\LL^1_\textrm{loc}([0,T];\LL^\infty(\RR^d))$.
\end{proof}

\bigskip

\begin{remark}
It may be tempting, in the light of the results exposed in Section~\ref{SubsectionIntegrationConservative}, to try and develop an $\LL^\infty$ theory of differential equations driven by more general rough drivers ${\bfA}_{ts}$ than those associated with the data of some vector fields $V_1,\dots, V_\ell$ and a weak geometric H\"older rough path over $\RR^\ell$, as in the transport equation \eqref{EqTransportEq}. With a view towards the classical theory of stochastic flows, as developed by Le Jan-Watanabe, Kunita and many others, one may try, as a first step, to work with rough drivers whose first level are obtained as typical trajectories of semimartingale velocity fields. It is shown in \cite{BailleulRiedel} that such random fields can be lifted into some objects very similar to rough drivers, under some mild regularity conditions on the semimartingale, and that the use of the approximate flow machinery introduced in \cite{BailleulRMI} leads to some well--posedness result for some dual evolution equation 
$$
df_t(\varphi) = f_t({\bf A}(dt)\varphi) .   
$$ 
\end{remark}

\end{document}